\documentclass[preprint,10.5pt]{elsarticle}
\pdfoutput=1

\usepackage[margin=1in]{geometry}

\usepackage{tikz}
\usepackage{array}
\usepackage{amsmath}
\usepackage{amssymb}
\usepackage{booktabs}
\usepackage{subcaption}
\usepackage{soul}
\usepackage{xcolor}
\usepackage{subcaption}
\usepackage{multirow}
\usepackage[export]{adjustbox}
\usepackage{ulem}
\usepackage{algorithm}
\usepackage{algpseudocode}
\usepackage{todonotes}
\usepackage{nicematrix}
\usepackage{cancel}
\usepackage{svg}

\def\NoNumber#1{{\def\alglinenumber##1{}\State #1}\addtocounter{ALG@line}{-1}}

\begin{document}

\begin{frontmatter}

\title{Online learning in bifurcating dynamic systems \\ via SINDy and Kalman filtering}

\author[1]{Luca Rosafalco\corref{cor1}}
\ead{luca.rosafalco@polimi.it}
\author[1,2]{Paolo Conti}
\ead{paolo.conti@polimi.it}
\author[2]{Andrea Manzoni}
\ead{andrea1.manzoni@polimi.it}
\author[1]{Stefano Mariani}
\ead{stefano.mariani@polimi.it}
\author[1]{Attilio Frangi}
\ead{attilio.frangi@polimi.it}
\address[1]{Department of Civil and Environmental Engineering,
Politecnico di Milano \\ Piazza L. da Vinci 32,
20133 - Milano (Italy)}
\address[2]{MOX, Department of Mathematics,
Politecnico di Milano \\ Piazza L. da Vinci 32,
20133 - Milano (Italy)}

\cortext[cor1]{Corresponding author}

\begin{abstract}
We propose the use of the Extended Kalman Filter (EKF) for online data assimilation and update of a dynamic model, preliminary identified through the Sparse Identification of Nonlinear Dynamics (SINDy). This data-driven technique may avoid biases due to incorrect modelling assumptions and exploits SINDy to approximate the system dynamics leveraging a predefined library of functions, where active terms are selected and weighted by a sparse set of coefficients. This results in a physically-sound and interpretable dynamic model allowing to reduce epistemic uncertainty often affecting machine learning approaches. Treating the SINDy model coefficients as random variables, we propose to update them while acquiring (possibly noisy) system measurements, thus enabling the online identification of time-varying systems. These changes can stem from, e.g., varying operational conditions or unforeseen events. The EKF performs model adaptation through joint state-parameters estimation, with the Jacobian matrices required to computed the model sensitivity inexpensively evaluated from the SINDy model formulation. The effectiveness of this approach is demonstrated through three case studies: (i) a Lokta-Volterra model in which all parameters simultaneously evolve during the observation period; (ii) a Selkov model where the system undergoes a bifurcation not seen during the SINDy training; (iii) a MEMS arch exhibiting a 1:2 internal resonance. The ability of EKF of recovering inactivated functional terms from the SINDy library, or discarding unnecessary contribution, is also highlighted. Based on the presented applications, this method shows strong promise for handling time-varying nonlinear dynamic systems possibly experiencing bifurcating behaviours.
\end{abstract}

\begin{keyword}
Time-varying nonlinear dynamic systems, Extended Kalman Filter (EKF), Sparse Identification of Nonlinear Dynamics (SINDy), Bifurcating dynamic systems.
\end{keyword}

\end{frontmatter}

\section{Introduction}
\label{sec:introduction}
The increasing availability of data has driven significant advancements across various scientific and engineering domains, including structural health monitoring in civil engineering \cite{book:Farrar13}, aerodynamic control \cite{art:Fonzi20}, and modeling the structural mechanics of Micro Electro-Mechanical Systems (MEMS) \cite{book:Corigliano18}. In all these fields,  data-driven approaches have emerged as new paradigms to analyse, forecast, and control complex physical systems, setting the foundation of predictive Digital Twins (DTs) \cite{art:Wagg20}, nowadays enabling, e.g., predictive maintenance for civil infrastructures \cite{art:Matt24} or unmanned aerial vehicles \cite{art:Kapteyn21}.

The construction of predictive DTs, unlocking the possibility to operate online on a dynamical system, requires to assimilate incoming information on-the-fly. This is possible thanks to filtering techniques like the Kalman Filter (KF) \cite{art:Kalman60}, that provides the optimal way to assimilate information in a linear setting \cite{book:Simon06_5}. Nonlinear extensions of the filter, such as the Extended Kalman Filter (EKF) or the unscented Kalman filter can be used for nonlinear identification problems \cite{proc:Wan00,art:Mariani07}. A broader application of filtering techniques for DTs is hampered by challenges in filter tuning and errors in formulating the dynamic model employed in the predictor stage of the filter \cite{art:RUENG24,art:Novoa24}. While advanced filter tuning approaches have been recently proposed in the literature \cite{art:Song20,art:Gres25,art:Bilgin25}, here we focus on data-driven non-intrusive surrogates to face the improper modelling issue suffered by purely physics-based models \cite{art:Cuomo24,art:Impraimakis24}.

Recent studies completely replaced physics-based models with Machine Learning (ML)-based models. For example, generative adversarial networks have been demonstrated successful in structural health monitoring \cite{art:SDEE23} and damage detection for nonlinear dynamic system \cite{art:Joseph24} in simulated environments. In a similar direction, the combination of ML-based models with KFs was recently explored. For instance, in \cite{proc:Coskun17}, long short-term memory networks were used in the predictor phase of the filter, resulting in the automatic estimation of the process and measurement noise. The work of \cite{art:Liu24} demonstrated impressive estimation capabilities in nonlinear structural dynamics by employing dynamical variational autoencoders whose inference task was performed by the EKF.

On the other hand, training Neural Networks (NNs) require large datasets and struggle with epistemic uncertainties \cite{art:Olivier21}, which refer to the inability of modelling a system response in scenarios where data are unavailable \cite{art:Cicirello24}. This is the case of civil structure responses in potential damage conditions \cite{chp:CSAI22}, where only baseline (undamaged) data are available. Such uncertainties, combined with the lack of interpretability and of physical consistency, cause significant limitations to the practical application of ML models in real engineering contexts, potentially hindering their use in the experimental setting \cite{art:Seventekidis23}. Moreover, pre-trained NNs generally cannot be updated online, potentially leading to considerable prediction errors in nonlinear dynamics. This is particularly problematic when highly nonlinear processes are involved, as dynamic contributions that were negligible during the initial NN training may emerge. For example, variations of input parameters can result in the development of vorticity in fluids \cite{conti2023reduced} or resonance in mechanical structures \cite{art:Opreni23}, possibly leading to bifurcation phenomena that completely change the system dynamics. In a DT perspective, 
detecting at which configurations these regime transitions occur is essential, as these behaviours can significantly affect the operation of the system. We highlight physical consistency as a required condition to make ML models more robust.


In this work, we present a data-driven method for online learning and adaptation of a dynamical model, using assimilated data from the physical system. Specifically, we adopt the Sparse Identification of Nonlinear Dynamics (SINDy) technique \cite{art:Brunton16} to construct the predictive model required in the filtering prediction stage; data assimilation is then carried out using an EKF. SINDy approximates the system dynamics from available data by constructing a nonlinear dynamic model that is parametrised by a sparse set of coefficients, determining the active terms from a predefined library of functions. The construction of robust, physically consistent, and interpretable models is thus enabled by SINDy thanks to its capacity to adopt a parsimonious description of system dynamics, consistent with the assumption that system models are governed by few important terms if an appropriate reference system is adopted, even when discontinuous nonlinearities are considered \cite{art:Lathourakis24}.

Differently from our previous work \cite{art:CMAME24}, here we do not only update system variables to track the evolution of states and parameters; rather, we evolve the predictive model itself by dealing with SINDy coefficients as learnable random variables. The outcome is an explainable data-driven model that can evolve online adapting to different scenarios. This allows to enrich a dynamic description of the physical phenomenon at hand as new data are assimilated, and to correct or even dropping terms that are incorrect or irrelevant. The capacity to assimilate data and modify online the SINDy model sets apart our work from previous applications of SINDy, for which training is typically done once in a preliminary offline phase.

With respect to \cite{art:CMAME24}, we no longer require a large amount of trajectories of the system for different parameter configurations to construct, a priori, a predictive model which accounts for parametric dependency. Recovering this type of data is not always possible or it may result extremely expensive, although possible \cite{art:Dardeno24}. Instead we can learn and incorporate this parametric dependency online, thus significantly enlarging the applicability of the method.

Thanks to the convenient non-intrusive system parametrisation provided by SINDy, the obtained framework can model parametric bifurcations, which is a fundamental and desirable goal. For example, the possibility of generalising the system dynamics over bifurcation parameter values was recently commented in \cite{art:Li24}, by adapting a NN-based meta-model of the system dynamics. However, this approach required to preliminary collect trajectories under varying the properties of interest.

Other works combined the ideas underlying the SINDy technique with filtering. For instance, in \cite{art:Gotte23,art:Pal24}, a sparsity promoting UKF algorithm was proposed for nonlinear dynamic identification, employing pseudo-observations as described in \cite{art:Julier07}, with the advantage of eliminating the need of any preliminary offline procedure; however, omitting SINDy for parameter initialisation makes the identification task considerably more challenging. In \cite{art:Stevens24}, the combination with a KF was proposed to mitigate the impact of noise on the identification of the SINDy model.

Similarly to what we propose, in \cite{proc:Wang22}, SINDy parameters were evolved, but the use of a linear version of the KF excluded the possibility of joint estimation of system dynamics and model parameters, thus precluding nonlinear dynamics applications for which being able to closely track the evolution of the system is essential. By adopting the EKF, we have been instead able to dramatically widen the range of applications of the procedure targeting a collection of interesting and challenging problems. In particular, we treat a case in which all the parameters of a dynamic model vary simultaneously during the observation window; a further case in which the variation of one parameter causes the system to pass through an Hopf bifurcation; a final situation related with a real MEMS device featuring a 1:2 internal resonance and Neimark-Sacker bifurcations \cite{art:Giorgio21}.

\begin{figure}[t!]
\centering
\includegraphics[width=145mm]{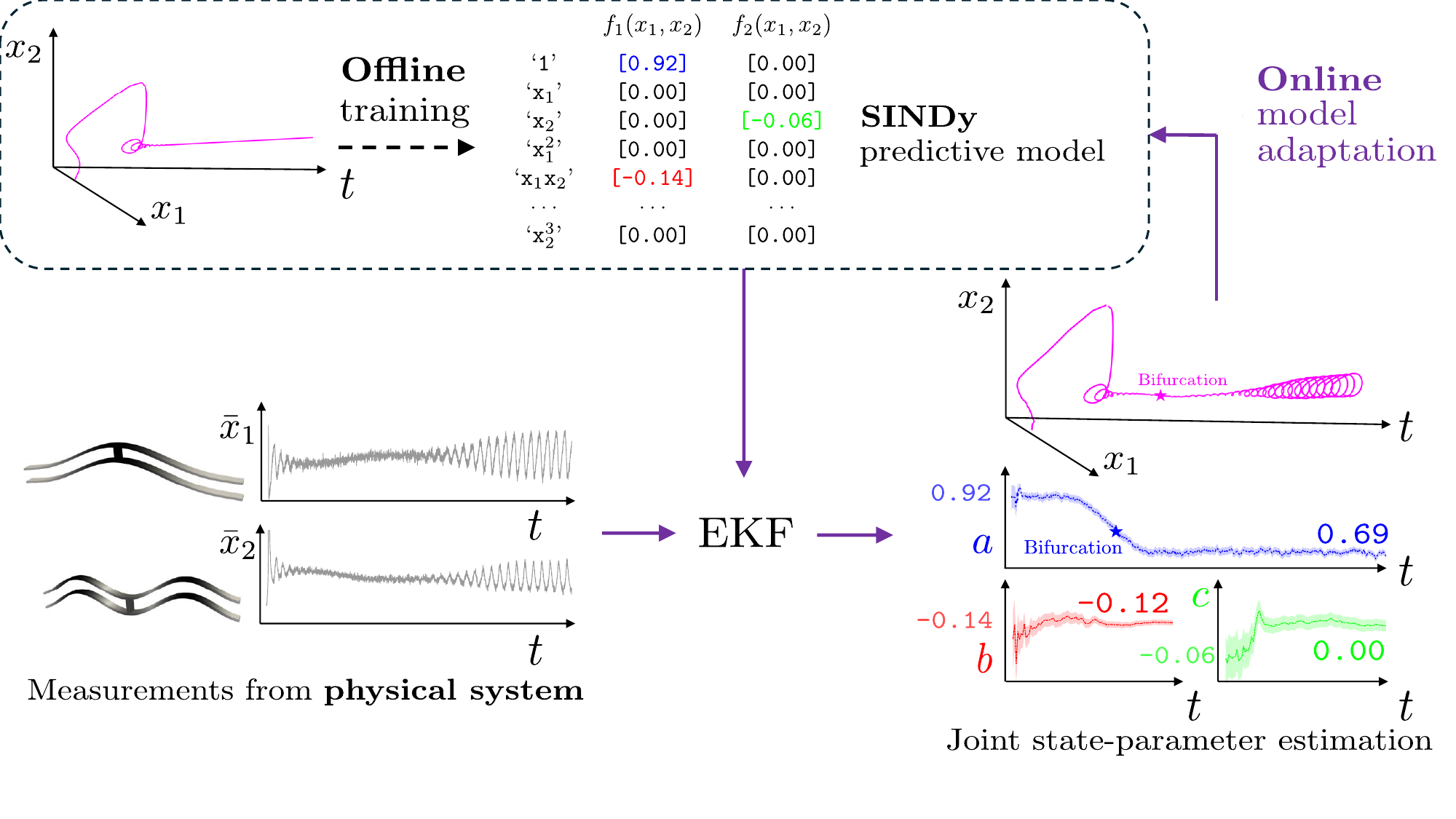}
{\caption{\footnotesize 
Online model adaptation combining SINDy and EKF. In the offline stage, a SINDy model is trained. In the online stage, this model is used in the predictor phase of the EKF to advance the system state over time. The EKF corrector phase then updates both the states and parameters of the model by assimilating measurements acquired on the physical system. The procedure can potentially recover inactivated functional terms from the SINDy library or, conversely, suppress unnecessary contributions.}\label{fig:graph_abs}}
\end{figure}

The remainder of the paper is organised as follows. In Sec. \ref{sec:methodology} the proposed procedure is discussed, focusing on EKF and SINDy equations, and highlighting how to combine them. A schematic representation of the proposed methodology, distinguishing between the offline SINDy training and its online update, is reported in Fig. \ref{fig:graph_abs}. Case studies are discussed in Sec. \ref{sec:results}, by addressing a time-varying Lokta-Volterra model (\ref{sec:LV}), a Selkov model passing an Hopf bifurcations (Sec. \ref{sec:Selkov}) and a real MEMS device (Sec. \ref{sec:MEMSarch}). Further comments and future developments are reported in Sec. \ref{sec:conclusions}.

\section{Methodology}
\label{sec:methodology}

We rely on a state-space representation of the dynamical system of interest. Such representation allows for the time evolution of the state variables describing the system dynamics, collected in the vector $\mathbf{x}(t)\in\mathbb{R}^n$; data assimilation exploits instead incoming observations of the system $\mathbf{y}(t_j)\in\mathbb{R}^{o}$, where $o$ are the number of observations at $t_j$, and $t_j$ with $j=1,\ldots,T$ are the discrete observation times. Observations and state components are treated as random variables, to account for modelling and measurement uncertainties. The following system of equations terms capture system evolution and observables:
\begin{subequations}
    \begin{align}
    \dot{\mathbf{x}}(t) = \frac{\text{d}}{\text{d} t}\mathbf{x}(t) = \mathbf{f}(\mathbf{x}(t),t) + \mathbf{q}_x, \quad t &\in [t_0,t_T], \label{eq:dynamicSystem} \\
    \mathbf{y}(t_j) = \mathbf{h}(\mathbf{x}(t_j)) + \mathbf{r}, \quad t_j &\in \lbrace t_0,t_1,\ldots,t_T \rbrace, \label{eq:obsEqDiscrete}
    \end{align}
\end{subequations}
where: $t$ is the time coordinate; $\mathbf{f}:\mathbb{R}^n\times \mathbb{R}\rightarrow \mathbb{R}^n$ is the function describing the system time-varying dynamics; $\mathbf{q}_x$ is the process noise vector; $\mathbf{h}:\mathbb{R}^n\rightarrow \mathbb{R}^o$ is a possibly nonlinear function mapping the state variables into the incoming measurements $\mathbf{y}(t_j)$. Both $\mathbf{q}_x$ and $\mathbf{r}$ are assumed to be distributed as white, zero mean stochastic processes, that is $\mathbf{q}_x\sim \mathcal{N}(\mathbf{0},\mathbf{Q}_x)$ with diagonal covariance matrix $\mathbf{Q}_x\in\mathbb{R}^{n\times n}$, and $\mathbf{r}\sim \mathcal{N}(\mathbf{0},\mathbf{R})$ with diagonal covariance matrix $\mathbf{R}\in\mathbb{R}^{o\times o}$, respectively. Eqs. \eqref{eq:dynamicSystem} and \eqref{eq:obsEqDiscrete} can be easily modified to include the dependence on a forcing term, as shown in \cite{art:CMAME24}. It is worth noting that the continuous time dependence in Eq. \eqref{eq:dynamicSystem} enables to directly substitute $\mathbf{f}$ with SINDy, while the discrete time dependence of the observables is consistent with digital recording acquisition.

A prediction-correction scheme is adopted: the quantities evolved through Eq.~\eqref{eq:dynamicSystem}, denoted with the superscript ``$-$", are the output of the prediction phase; these outputs will be then updated in the correction phase. The resulting quantities will be denoted with the superscript ``$+$". According to the extended version of the filter (EKF), the mean $\hat{\mathbf{x}}=\mathbb{E}[\mathbf{x}]$ and covariance $\mathbf{P}=\mathbb{E}[(\mathbf{x}-\hat{\mathbf{x}})(\mathbf{x}-\hat{\mathbf{x}})^\top]$ of the state variables $\mathbf{x}$ are propagated in time.

During the prediction stage, the values of $\hat{\mathbf{x}}^+(t_{j\!-\!1})$ and $\mathbf{P}^+(t_{j\!-\!1})$ at $t_{j\!-\!1}$ are evolved to obtain $\hat{\mathbf{x}}^-(t_j)$ and $\mathbf{P}^-(t_j)$ at $t_j$, respectively. To integrate the continuous relation in Eq.~\eqref{eq:dynamicSystem}, different methods can be utilised. Due to its high accuracy, we have employed the fourth-order Runge Kutta integrator. More details on time integration can be found in Algorithm \ref{al:onPhase}. To compute $\mathbf{P}^+(t_j)$, the Jacobian $\mathbf{F}^{+}(t_{j\!-\!1})$ of $\mathbf{f}$ has to be computed according to:
\begin{equation}
\mathbf{F}^{+}(t_{j\!-\!1}) =\frac{\partial \mathbf{f}}{\partial \mathbf{x}^{\top}}\bigg|_{\hat{\mathbf{x}}^+_{(t_{j\!-\!1})}}.
\label{eq:Jacobian_def}
\end{equation}

After getting $\hat{\mathbf{x}}^+(t_{j\!-\!1})$ and $\mathbf{P}^+(t_{j\!-\!1})$ in the prediction stage, their estimates are updated by assimilating the information from incoming measurements $\mathbf{y}(t_j)$ to get of $\hat{\mathbf{x}}^+(t_j)$ and $\mathbf{P}^+(t_j)$ at $t_j$. This is achieved by way of:
\begin{subequations}
    \begin{align}
        \mathbf{G}(t_j) = \mathbf{P}^-(t_j) {\mathbf{H}^-(t_j)}^\top (\mathbf{H}^-(t_j) \mathbf{P}^-(t_j) {\mathbf{H}^-(t_j)}^\top \!\!+ \mathbf{R})^{-1}, \\
        \hat{\mathbf{x}}^+(t_j) = \hat{\mathbf{x}}^-(t_j) + \mathbf{G}(t_j) (\mathbf{y}(t_j) - \mathbf{h}(\hat{\mathbf{x}}^-(t_j),t_j)),
        \label{eq:statePosterior}\\
        \mathbf{P}^+\!(t_j)\! =\! \left(\mathbf{I}\!-\!\mathbf{G}(t_j)\mathbf{H}^-(t_j)\right)\mathbf{P}^-(t_j)\left(\mathbf{I}\!-\!\mathbf{G}(t_j) \mathbf{H}^-(t_j)\right)^\top\!\! \!+\! \mathbf{G}(t_j) \mathbf{R} \mathbf{G}(t_j)^\top,
        \label{eq:covPosterior}
    \end{align}
\end{subequations}
where: $\mathbf{G}(t_j)$ is the Kalman gain at time $t_j$; $\mathbf{H}$ is the Jacobian of $\mathbf{h}$, computed as:
\begin{equation}
    \mathbf{H}^-(t_j) =\frac{\partial \mathbf{h}}{\partial \mathbf{x}^{\top}}\bigg|_{\hat{\mathbf{x}}^-(t_j)}.
    \label{eq:obsEqsJac}
\end{equation}

In this work, we assume $\mathbf{h}$ to be an operator selecting which state variable must be compared with the incoming observations. Alternatively, the same identification strategy here adopted for $\mathbf{f}$ can be employed to identify not trivial $\mathbf{h}$, as already proposed in \cite{art:CMAME24}. 

To avoid modelling issues \cite{art:Novoa24}, we identify $\mathbf{f}$ directly from data. In particular, we consider the SINDy setup \cite{art:Brunton16}, where the function $\mathbf{f}$ is approximated by a linear combination of $p$ candidate functions included in a library $\boldsymbol{\Theta}(\mathbf{x}(t))=\left[\theta_1(\mathbf{x}(t)), \ldots, \theta_p(\mathbf{x}(t))\right]\in\mathbb{R}^p$ that describes the system dynamics: 
\begin{equation}
    \mathbf{f}(\mathbf{x}(t),t)\approx \mathbf{\Xi}^\top(t)\mathbf{\Theta}^\top(\mathbf{x}(t)),
    \label{eq:f_sindy_approx}
\end{equation}
where $\boldsymbol{\Xi}(t) =\left[ \boldsymbol{\xi}_1(t),\ldots, \boldsymbol{\xi}_n(t)\right] \in \mathbb{R}^{p\times n}$ is the matrix collecting the weighting coefficients $\boldsymbol{\xi}_i(t)\in \mathbb{R}^{p}$.  Any type of functions can be included in $\boldsymbol{\Theta}$, with polynomials often serving as a suitable library for problems in computational mechanics \cite{art:Champion19}. It is also possible to assume that the functions in $\boldsymbol{\Theta}$ have an explicit time dependency, thus getting $\mathbf{f}(\mathbf{x}(t),t)\approx \mathbf{\Xi}^\top(t) \mathbf{\Theta}^\top(\mathbf{x}(t),t)$, for example to treat non-autonomous (externally forced) systems \cite{art:CMAME24}.

Our goal is to estimate the values of the unknown coefficients $\boldsymbol{\Xi}$ and their time evolution from the assimilated data. It is important to note that this framework allows non-zero coefficients in $\boldsymbol{\Xi}$ to vanish or, conversely, zero-valued coefficients to increase over time. This flexibility enables potential changes in the composition of the functional terms that contribute to the system dynamics.

First, an initialisation of the coefficients $\boldsymbol{\Xi}$ is performed offline using static learning procedures on a limited training set of snapshots of the state vector trajectories. The training data are collected in the following matrix:
\begin{equation}
    \mathbf{X} = \begin{bmatrix}        
    \mathbf{x}^{\top}(t_1) \\
    \mathbf{x}^{\top}(t_2) \\
    \vdots \\
    \mathbf{x}^{\top}(t_T)
    \end{bmatrix} = 
    \begin{bmatrix}
        x_1(t_1) & x_2(t_1) & \cdots & x_n(t_1) \\
        x_1(t_2) & x_2(t_2) & \cdots & x_n(t_2) \\
        \vdots & \vdots & \ddots & \vdots \\
        x_1(t_T) & x_2(t_T) & \cdots & x_n(t_T)
    \end{bmatrix},
    \label{eq:snapshotMatrix}
\end{equation}
where $x_i(t_j)$ is the i-th entry of $\mathbf{x}$ at the $j$-th time instant. Similarly, a matrix $\dot{\mathbf{X}}\in\mathbb{R}^{T\times n}$ is constructed by collecting the derivatives $\dot{\mathbf{x}}$. Thus, the functions of the library are applied to the rows of $\mathbf{X}$ obtaining $\boldsymbol{\Theta}(\mathbf{X})\in \mathbb{R}^{T\times  p}$ (see \cite{art:Brunton16} for further details) and Eq.~\eqref{eq:f_sindy_approx} translates at the levels of the training data as:
\begin{equation}
    \dot{\mathbf{X}} \approx\boldsymbol{\Theta}(\mathbf{X}) \boldsymbol{\Xi}.
    \label{eq:sindy}
\end{equation}
Matrix $\boldsymbol{\Xi}$ is obtained by solving the following sparsity promoted least square regression problem:
\begin{equation}
\boldsymbol{\Xi}=\underset{\boldsymbol{\Xi}'}{\arg\min}\Vert \dot{\mathbf{X}} - \boldsymbol{\Theta}(\mathbf{X})\boldsymbol{\Xi}' \Vert_2 + \mathcal{L}(\boldsymbol{\Xi}'), 
\label{eq:weightSINDy}
\end{equation}
where $\mathcal{L}(\boldsymbol{\Xi}')$ is a sparsity promoting regulariser that could enforce, for example, a LASSO regression \cite{hastie2015statistical} or a Sequential Thresholded Least SQuares (STLSQ) algorithm \cite{art:Brunton16}. Notice that this intialization represents a static estimate, representing the best fitting constant coefficients over the set of training trajectories.

Next, we are interested to adapt $\mathbf{f}$ to assimilate acquired data, in order to account for time-varying dynamics. 
To do so, we augment the state vector to include a subset of the SINDy weighting coefficients $\tilde{\boldsymbol{\xi}}$, thus obtaining $\boldsymbol{\varkappa}=[\mathbf{x},\tilde{\boldsymbol{\xi}}]^{\top}$. This means that the SINDy coefficients of this subset are treated as random variables. In other words, $\tilde{\boldsymbol{\xi}}$ collects the entries $\boldsymbol{\Xi}$ that may evolve in the system, while the remaining coefficients are grouped into $\boldsymbol{\xi}_0$. Since many coefficients are usually set to zero through Eq.~\eqref{eq:weightSINDy}, we expect that many of them will remain unchanged. On the contrary, coefficients set to zero by SINDy may be included in $\tilde{\boldsymbol{\xi}}$ to track potential emerging dynamic behaviours. Coefficients $\tilde{\boldsymbol{\xi}}$ and $\boldsymbol{\xi}_0$ are collected into $\tilde{\mathbf{\Xi}}$ and $\mathbf{\Xi}_0\in\mathbb{R}^{p\times n}$, respectively, by multiplying $\mathbf{\Xi}$ through a Boolean matrix $\mathbf{B}\in\mathbb{R}^{p\times n}$ elementwise, that is by multiplying corresponding elements of $\mathbf{B}$ and $\mathbf{\Xi}$:
\begin{subequations}
    \begin{align}
        \mathbf{\Xi}(t) = \tilde{\mathbf{\Xi}}(t) + \mathbf{\Xi}_0, \quad \text{with}\\
        \tilde{\mathbf{\Xi}}(t) = \mathbf{B}\odot\mathbf{\Xi}(t),\\
        \mathbf{\Xi}_0 = (\mathbf{1}-\mathbf{B})\odot\mathbf{\Xi}(t),
    \end{align}
\end{subequations}
where $\odot$ denotes the elementwise product operator and $\mathbf{1}\in\mathbb{R}^{p\times n}$ is a matrix whose entries are set to $1$.

Consistently, the introduced notation allows the dynamic evolution of $\mathbf{x}$ to be expressed as:
\begin{equation}
    \mathbf{f}(\mathbf{x}(t))\approx \mathbf{\Xi}(t)^\top\mathbf{\Theta}^\top(\mathbf{x}(t),t)=(\tilde{\mathbf{\Xi}}(t)+\mathbf{\Xi}_0)^\top\mathbf{\Theta}^\top(\mathbf{x}(t)).
\end{equation}
However, defining the augmented state $\boldsymbol{\varkappa}$ requires describing the dynamic evolution of $\tilde{\boldsymbol{\xi}}$ as well. We assume that this evolution is driven by a white, zero mean stochastic process noise $\mathbf{q}_{\hat{\xi}}\sim \mathcal{N}(\mathbf{0},\mathbf{Q}_{\hat{\xi}})$ with diagonal covariance matrix $\mathbf{Q}_{\hat{\xi}}\in\mathbb{R}^{\varrho\times \varrho}$. We thus introduce a unique process noise $\mathbf{q}=[\mathbf{q}_x,\mathbf{q}_{\hat{\xi}}]^{\top}$, and a unique covariance matrix $\mathbf{Q}$ collecting the entries of $\mathbf{Q}_x$ and $\mathbf{Q}_{\hat{\xi}}$. This means that the values of $\tilde{\boldsymbol{\xi}}$ can be only modified in the correction stage on the basis of the mismatch between the incoming $\mathbf{y}(t_j)$ and the predicted $\mathbf{h}(\hat{\mathbf{x}}^-(t_j),t_j)$.

Ultimately, using SINDy to estimate $\mathbf{f}$ enables to compute the Jacobian matrix $\mathbf{F}^{+}(t_{j\!-\!1})$ in a straightforward way as detailed in the following:
\begin{equation}
\mathbf{F}^{+}(t_{j\!-\!1})\! = \!\begin{bmatrix}
    \frac{\partial \mathbf{f}}{\partial \mathbf{x}^{\top}}\bigg|_{\hat{\boldsymbol{\varkappa}}^{+}(t_{j\!-\!1})} & \frac{\partial \mathbf{f}}{\partial \tilde{\boldsymbol{\xi}}^{\top}}\bigg|_{\hat{\boldsymbol{\varkappa}}^{+}(t_{j\!-\!1})} \\
    \mathbf{0}   & \mathbf{0}
\end{bmatrix}\!=\!\begin{bmatrix}
(\tilde{\boldsymbol{\Xi}}^+(t_{j\!-\!1}) + \boldsymbol{\Xi}_0) ^{\top}\frac{\partial\boldsymbol{\Theta}^{\top}(\mathbf{x})}{\partial \mathbf{x}^{\top}}\bigg|_{\hat{\boldsymbol{\varkappa}}^{+}(t_{j\!-\!1})} & \frac{\partial (\tilde{\boldsymbol{\Xi}}(t_{j\!-\!1}))^{\top}}{\partial \tilde{\boldsymbol{\xi}}}\tilde{\boldsymbol{\Theta}}^{\top}(\mathbf{x})\bigg|_{\hat{\boldsymbol{\varkappa}}^{+}(t_{j\!-\!1})}  \\
\mathbf{0} & \mathbf{0}
\end{bmatrix},
\label{eq:Jacobian_def_augmented}
\end{equation}
where it can be noted how the the derivatives with respect to $\tilde{\boldsymbol{\xi}}$ directly enters in the Jacobian definition. Implementation details of how the online adaption of the SINDy model is performed via the EKF by assimilating incoming measurements is reported in Algorithm \ref{al:onPhase}.

\begin{algorithm}[t!]
\hspace*{\algorithmicindent} \textbf{Input}: SINDy model $\mathbf{f}(\mathbf{x}(t))\approx \mathbf{\Xi}(t)^\top\mathbf{\Theta}^\top(\mathbf{x}(t))=(\tilde{\mathbf{\Xi}}(t)+\mathbf{\Xi}_0)^\top\mathbf{\Theta}^\top(\mathbf{x}(t))$; selecting operator $\mathbf{h}(\mathbf{x}(t_j),t_j)$; sequential measurements $\mathbf{y}(t_{1}),\ldots,\mathbf{y}(t_{T})$ \\
\hspace*{\algorithmicindent} \textbf{Output}: online update of the dynamic state variables $\mathbf{x}$ and of the SINDy coefficients $\tilde{\boldsymbol{\xi}}$ through the estimates of $\hat{\boldsymbol{\varkappa}}^+(t_1),\ldots,\hat{\boldsymbol{\varkappa}}^+(t_T)$, and $\hat{\mathbf{P}}^+(t_1),\ldots,\hat{\mathbf{P}}^+(t_T)$.\\
\hspace*{\algorithmicindent} \textbf{Set}: process and measurement noise covariance matrices $\mathbf{Q}$ and $\mathbf{R}$; the Boolean matrix $\mathbf{B}$.
\begin{algorithmic}[1]
\State $\mathbf{\Xi}(t)=\tilde{\mathbf{\Xi}}(t)+\mathbf{\Xi}_0$, with $\tilde{\mathbf{\Xi}}(t) = \mathbf{B}\odot\mathbf{\Xi}(t)$ and $\mathbf{\Xi}_0 = (\mathbf{1}-\mathbf{B})\odot\mathbf{\Xi}(t)$
\For{$t_j=t_1,\ldots,t_T$}
\[
\text{Predictor phase (4th order Runge Kutta integration)}
\]
\State $\hat{\mathbf{x}}^-(t_{j})=\hat{\mathbf{x}}^+(t_{j-1})+\frac{\Delta t}{6}(\mathcal{K}^{\mathbf{x}}_1 + 2\mathcal{K}^{\mathbf{x}}_2 + 2\mathcal{K}^{\mathbf{x}}_3+\mathcal{K}^{\mathbf{x}}_4)$ 
\NoNumber{}{\Comment{$\mathcal{K}^{\mathbf{x}}_1=(\tilde{\mathbf{\Xi}}^+(t_{j-1})+\mathbf{\Xi}_0)^\top\mathbf{\Theta}^\top(\hat{\mathbf{x}}^+(t_{j-1}))$}}
\NoNumber{}{\Comment{$\mathcal{K}^{\mathbf{x}}_2=(\tilde{\mathbf{\Xi}}^+(t_{j-1})+\mathbf{\Xi}_0)^\top\mathbf{\Theta}^\top(\hat{\mathbf{x}}^+(t_{j-1})\!+\!\Delta t \frac{\mathcal{K}^{\mathbf{x}}_1}{2})$}}
\NoNumber{}{\Comment{$\mathcal{K}^{\mathbf{x}}_3=(\tilde{\mathbf{\Xi}}^+(t_{j-1})+\mathbf{\Xi}_0)^\top\mathbf{\Theta}^\top(\hat{\mathbf{x}}^+(t_{j-1})\!+\!\Delta t \frac{\mathcal{K}^{\mathbf{x}}_2}{2})$}}
\NoNumber{}{\Comment{$\mathcal{K}^{\mathbf{x}}_4=(\tilde{\mathbf{\Xi}}^+(t_{j-1})+\mathbf{\Xi}_0)^\top\mathbf{\Theta}^\top(\hat{\mathbf{x}}^+(t_{j-1})\!+\!\Delta t \mathcal{K}^{\mathbf{x}}_3)$}}
\State $\hat{\tilde{\boldsymbol{\xi}}}^-(t_{j})=\hat{\tilde{\boldsymbol{\xi}}}^+(t_{j-1})$
\State $\mathbf{P}^-(t_{j})=\mathbf{P}^+(t_{j\!-\!1})+\frac{\Delta t}{6}(\mathcal{K}^{\mathbf{P}}_1 + 2\mathcal{K}^{\mathbf{P}}_2 + 2\mathcal{K}^{\mathbf{P}}_3+\mathcal{K}^{\mathbf{P}}_4)$
\NoNumber{}{\Comment{$\mathcal{K}^{\mathbf{P}}_1= \mathbf{F}^+(t_{j\!-\!1})\mathbf{P}^+(t_{j\!-\!1}) + \mathbf{P}^+(t_{j\!-\!1}) \mathbf{F}^{+^\top}\!(t_{j\!-\!1}) + \mathbf{Q}_x$}}
\NoNumber{}{\Comment{$\mathcal{K}^{\mathbf{P}}_2= \mathbf{F}^+(t_{j\!-\!1}+\frac{\Delta t}{2})\left(\mathbf{P}^+(t_{j\!-\!1})\!+\!\frac{\mathcal{K}^{\mathbf{P}}_1}{2}\right)\! + \!\left(\mathbf{P}^+(t_{j\!-\!1}) \!+\!\frac{\mathcal{K}^{\mathbf{P}}_1}{2}\right)\mathbf{F}^{+^\top}\!(t_{j\!-\!1}\!+\!\frac{\Delta t}{2})\! +\! \mathbf{Q}_x$}}
\NoNumber{}{\Comment{$\mathcal{K}^{\mathbf{P}}_3= \mathbf{F}^+(t_{j\!-\!1}\!+\!\frac{\Delta t}{2})\left(\mathbf{P}^+(t_{j\!-\!1})\!+\!\frac{\mathcal{K}^{\mathbf{P}}_2}{2}\right)\! + \!\left(\mathbf{P}^+(t_{j\!-\!1}) \!+\!\frac{\mathcal{K}^{\mathbf{P}}_2}{2}\right)\mathbf{F}^{+^\top}\!(t_{j\!-\!1}\!+\!\frac{\Delta t}{2}) \!+ \!\mathbf{Q}_x$}}
\NoNumber{}{\Comment{$\mathcal{K}^{\mathbf{P}}_4= \mathbf{F}^+(t_{j\!-\!1}\!+\!\frac{\Delta t}{2})\left(\mathbf{P}^+(t_{j\!-\!1})\!+\!\mathcal{K}^{\mathbf{P}}_3\right) \!+\! \left(\mathbf{P}^+(t_{j\!-\!1}) \!+\!\mathcal{K}^{\mathbf{P}}_3\right)\mathbf{F}^{+^\top}\!(t_{j\!-\!1}+\Delta t)\! +\! \mathbf{Q}_x$}}
\NoNumber{}{\Comment{$\mathbf{F}^{+}(t_{j\!-\!1}) = \begin{bmatrix}
(\tilde{\boldsymbol{\Xi}}^+(t_{j\!-\!1}) + \boldsymbol{\Xi}_0) ^{\top}\frac{\partial\boldsymbol{\Theta}^{\top}(\mathbf{x})}{\partial \mathbf{x}^{\top}}\bigg|_{\hat{\boldsymbol{\varkappa}}^{+}(t_{j\!-\!1})} & \frac{\partial (\tilde{\boldsymbol{\Xi}}^+(t_{j\!-\!1}))^{\top}}{\partial \tilde{\boldsymbol{\xi}}}\tilde{\boldsymbol{\Theta}}^{\top}(\mathbf{x})\bigg|_{\hat{\boldsymbol{\varkappa}}^{+}(t_{j\!-\!1})}  \\
\mathbf{0} & \mathbf{0}
\end{bmatrix}$}}
\[
\text{Corrector phase}
\]
\State $\mathbf{G}(t_j) = \mathbf{P}^-(t_j) {\mathbf{H}^-(t_j)}^\top (\mathbf{H}^-(t_j) \mathbf{P}^-(t_j) {\mathbf{H}^-(t_j)}^\top \!\!+ \mathbf{R})^{-1}$
\State $ \hat{\mathbf{x}}^+(t_j) = \hat{\mathbf{x}}^-(t_j) + \mathbf{G}(t_j) (\mathbf{y}(t_j) - \mathbf{h}(\hat{\mathbf{x}}^-(t_j),t_j))$
{\Comment{$\mathbf{H}^-(t_j) =\begin{bmatrix}
    \frac{\partial \mathbf{h}}{\partial \mathbf{x}^{\top}}\bigg|_{\hat{\mathbf{x}}^-(t_j)} & \mathbf{0}
\end{bmatrix}$}}
\State $\mathbf{P}^+(t_j)\! =\! \left(\mathbf{I}-\mathbf{G}(t_j)\mathbf{H}^-\!(t_j)\right)\mathbf{P}^-(t_j)\left(\mathbf{I}-\mathbf{G}(t_j) \mathbf{H}^-\!(t_j)\right)^\top \!+ \mathbf{G}(t_j) \mathbf{R} \mathbf{G}(t_j)^\top$
\EndFor
\end{algorithmic}
\caption{Online adaptation of SINDy via EKF}
\label{al:onPhase}
\end{algorithm}

\section{Numerical Results}
\label{sec:results}

The method outlined in Sec. \ref{sec:methodology} is now applied to three case studies. In Sec. \ref{sec:LV}, we examine a Lokta-Volterra model where all parameters evolve simultaneously, demonstrating our framework ability to construct and online update a SINDy model according to acquired measurements. In Sec. \ref{sec:Selkov}, we investigate a Selkov model to determine whether the adaptive SINDy model can track previously unseen dynamics after undergoing a bifurcation point controlled by an evolving parameter. All variables and parameters in these two case studies are treated as nondimensional. In Sec. \ref{sec:MEMSarch}, we assess the applicability of our approach to a pair of industry relevant engineering scenarios involving a MEMS device with a 1:2 internal resonance.

\subsection{Lokta Volterra model}
\label{sec:LV}

The Lokta-Volterra model is used to describe the dynamics of a biological system involving the interaction between two species, one acting as prey and the other as predator \cite{art:Lokta1920}. The system dynamics reads:
\begin{subequations}
    \begin{align}
        \dot{x}_1 = f_1(x_1,x_2)=  a x_1 + b x_1 x_2, \\
        \dot{x}_2 = f_2(x_1,x_2)=  c x_2 + d x_1 x_2 ,
        \label{eq:LV_equations}
    \end{align}
\end{subequations}
where $x_1$ is the prey population quantity, $x_2$ is the predator population quantity; $a\in\mathbb{R}^+$, $b\in\mathbb{R}^-$, $c\in\mathbb{R}^-$, $d\in\mathbb{R}^+$ are parameters characterising the interaction between the two species, and $\mathbf{f}=[f_1, f_2]^\top$.

We aim to evaluate whether the proposed procedure can track potential changes in the model parameters within the observation window $[0,T]$, with $T=150$, by monitoring both $x_1$ and $x_2$. To achieve this, we consider a function library $\boldsymbol{\Theta}$ consisting of polynomial terms up to the second order. The corresponding coefficients $\boldsymbol{\Xi}$ are preliminary identified offline using STLSQ over a set of system trajectories not affected by noise, computed by numerical integration for the following model parameters $a=1.0, b=-0.1, c=-1.5$ and $d=0.075$. The STLSQ algorithm is performed with a regularisation strength $\delta_r=0.05$, as default in the \texttt{PySINDy} package, and a threshold parameter $L=5\cdot10^{-4}$. Specifically, $\delta_r$ determines the strength of the $\ell$-2 regularisation promoting sparsity in the SINDy representation, while $L$ fixes the threshold below which the entries of $\boldsymbol{\Xi}$ are set to $0$. In the case of noisy data, noise tolerant versions of SINDy could be employed \cite{art:Fasel22,art:Paolo24}. As shown in \ref{sec:appendixA}, SINDy successfully reconstructs the original model as can be noted by comparing the coefficients determined by SINDy with the ones that weight the relevant functional terms in the reference model used to train SINDy.

After the offline initialisation of the model, we employ the EKF to online assimilate noisy observations $\mathbf{y} = [\bar{x}_1, \bar{x}_2]^{\top}$ of the system, letting the parameters evolve according to Algorithm \ref{al:onPhase}, and thereby consistently updating $\mathbf{f}$. This has involved incorporating $a$, $b$, $c$, $d$ into $\tilde{\boldsymbol{\xi}}$. The integration time step $\Delta t$, equal to the sampling time, has been set to $5.130\cdot 10^{-3}$. The acquired measurements have been characterised by a Signal to Noise Ratio (SNR) of $25$ with white noise generated numerically. To test the procedure in more challenging conditions, each parameter has been assigned a different evolution law. A sinusoidal variation, ranging from $0.8$ to $1.2$, has been considered for $a$, while a step change from $-0.1$ to $-0.09$ at $t=50$ have been imposed on $b$; $c$ has been instead kept constant to assess the filter ability to maintain its estimate unchanged; lastly, $d$ has been given a linear variation from $0.075$ to $0.085$. It is important to note that, by keeping unchanged the sign of the model coefficients, the stability portrait of the system remains unchanged, resulting in periodic trajectories which form closed orbits in the state-space surrounding the stable center $\mathbf{x}_{\text{eq}}= [\frac{d}{c},\frac{a}{b}]^\top$.

\begin{figure}[t!]
    \centering
    \subfloat[]{%
\includegraphics[width=145mm]{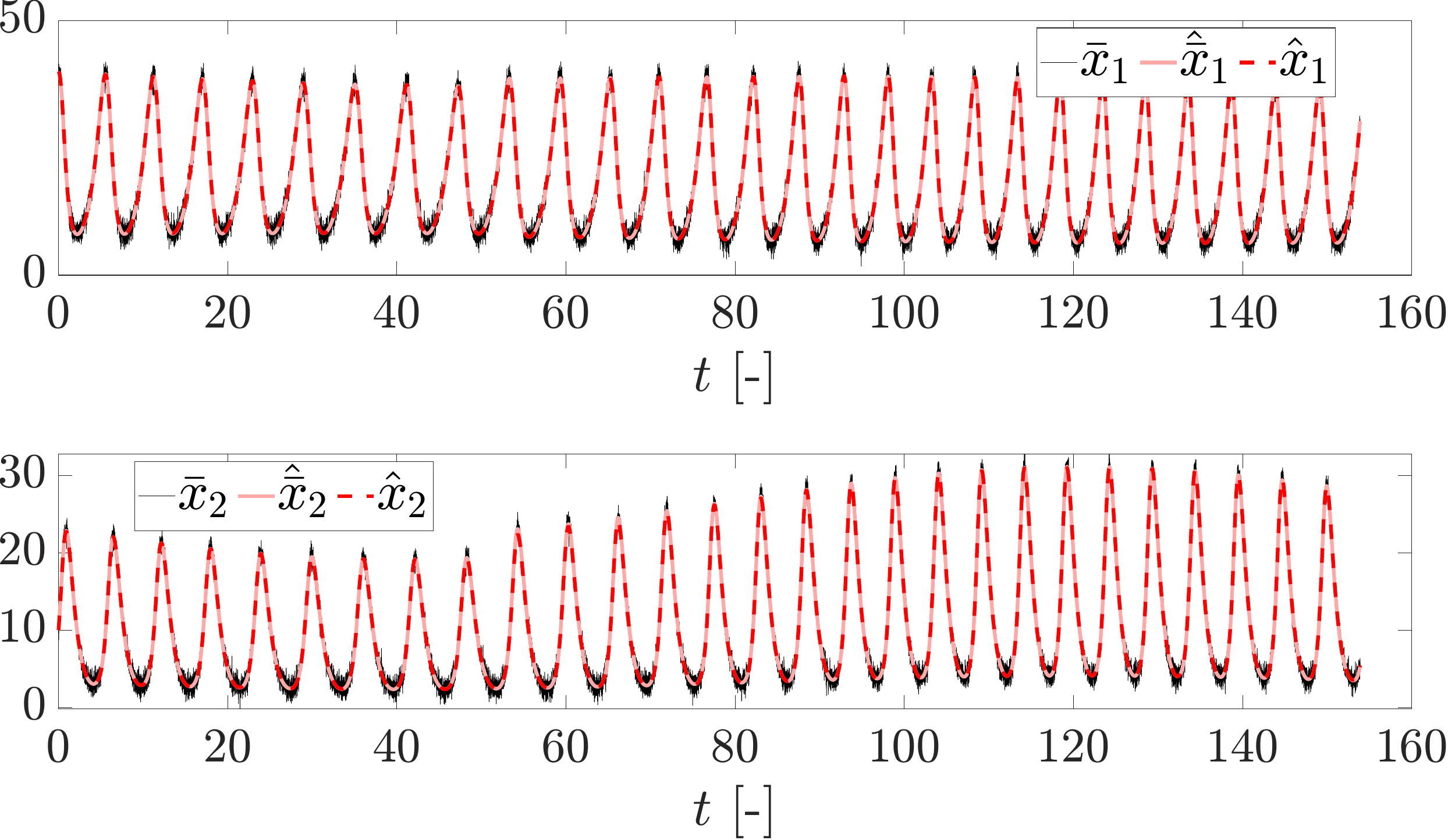}
\label{fig:LV_prey_predator}
} \\
 \subfloat[]{\includegraphics[width=145mm]{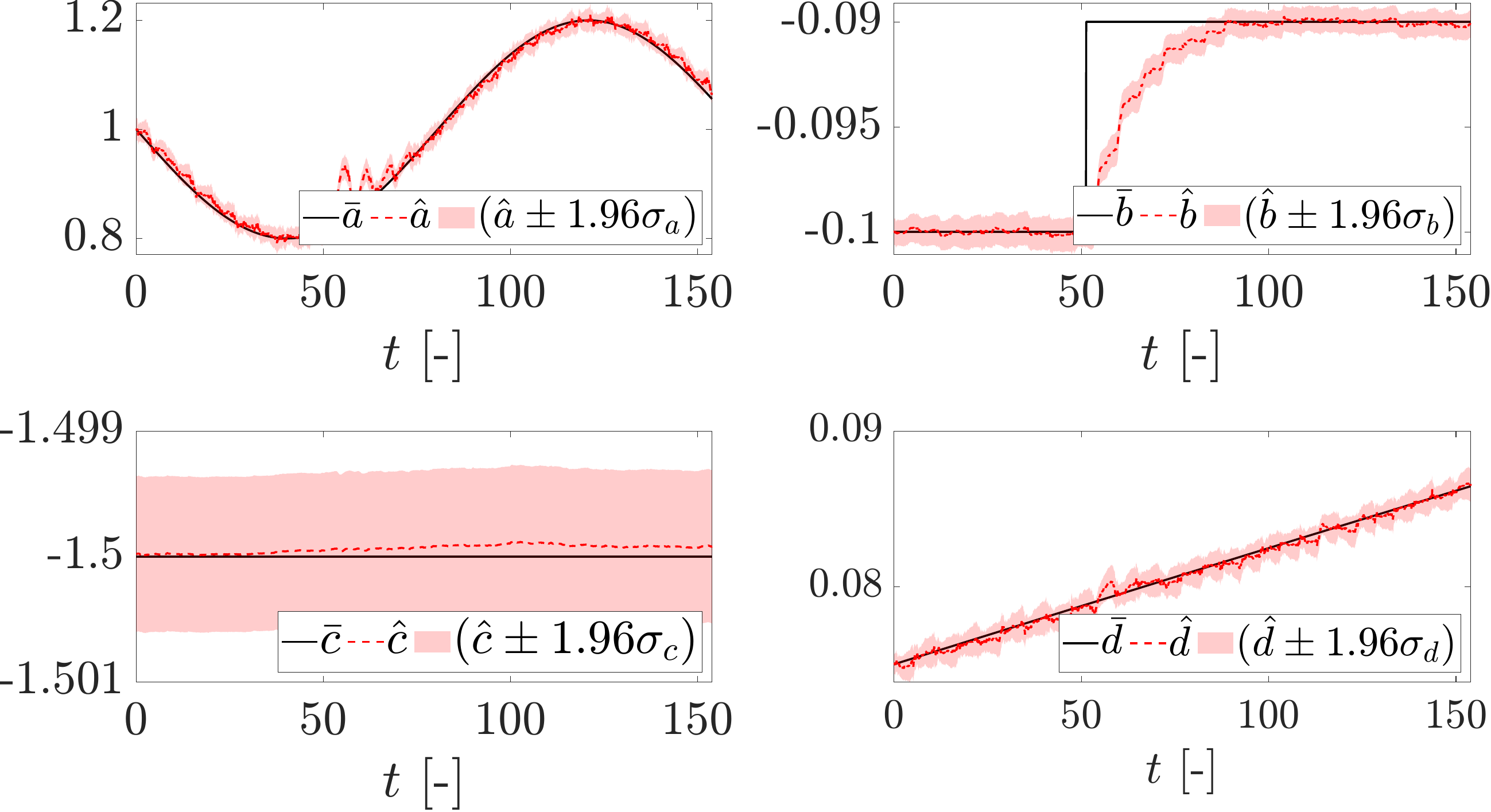}
\label{fig:LV_parameters}}
\caption{Lokta-Volterra model. In the top figure (a), the estimated states are plotted (red dotted line) against the assimilated noisy signals (black solid lines). Additionally, noise-free signals are reported in pink solid lines. In the bottom figure (b), the time-evolution of the estimated SINDy-model parameter are plotted (red dotted line) against the target time evolution (black solid line). The red shaded area represents 95\% confidence interval of the estimates, determined using the posterior covariance.}
\label{fig:LV}
\end{figure}

The results of the identification are collected in Fig. \ref{fig:LV} in terms both of tracking the dynamics of the system and of describing its parameter evolution. Target values are indicated with a bar, while noise-free version of the observations are labelled as $\hat{\bar{x}}_1$ and $\hat{\bar{x}}_2$. The procedure has allowed us to successfully track the system dynamics even under evolving conditions in the predictive model $\mathbf{f}$, which had not been encountered during SINDy training. All parameters have been accurately and rapidly identified, particularly in response of the step change in $b$, proving how model adaptability can be attained by assimilating incoming measurements through EKF.

\subsection{Selkov model}
\label{sec:Selkov}

The Selkov model describes glycolysis by the following system of two ordinary differential equations \cite{art:Selkov68}:
\begin{subequations}
\begin{align}
    \dot{x}_1 &= f_1(x_1,x_2) = \rho +\varsigma_1  x_1 -\vartheta_1 x_1  x_2^2 , \label{eq:selkov_1} \\
    \dot{x}_2 &= f_2(x_1,x_2) = \varsigma_2  x_1 + \psi x_2 + \vartheta_2 x_1 x_2^2 , \label{eq:selkov_2} 
\end{align}
\label{eq:selkov}
\end{subequations}
\noindent where $x_1$ and $x_2$ represent the concentration of two chemicals species in the glycolysis. Not all the coefficients that weight the functional terms usually vary when applying the Selkov model. Specifically, $\vartheta_1=-1$, $\vartheta_2=1$, $\psi=-1$ are typically fixed. Similarly,
the contribution of the cubic term $x_1 x^2_2$ is generally not parametrised, but it is included here with coefficients $\vartheta_1\in\mathbb{R}^-$ and $\vartheta_2\in\mathbb{R}^+$. Likewise, two separate parameters $\varsigma_1$ and $\varsigma_2$ have been used instead of the usual single coefficient $\varsigma$, which would normally require changing the sign of the $x_1$ term in Eq.~\eqref{eq:selkov_1}). This detailed parametrisation has been adopted for clarity, particularly to demonstrate whether SINDy, combined with the EKF, can estimate all these coefficients by observing $x_1$ and $x_2$.

Particular attention is given to $\rho$, as varying its value can alter the system dynamics. By setting $\varsigma = 0.1$, the system undergoes two Hopf bifurcations \cite{munoz2011introduction} at $\rho \simeq 0.42$ and $\rho \simeq 0.79$. Specifically, for $\rho < 0.42$ and $\rho > 0.79$ the system has a stable fixed point, while for $0.42 < \rho < 0.79$ the fixed point loses stability, giving rise to a stable limit cycle. In the following, we show how the EKF is capable of detecting and tracking the system through the Hopf bifurcation occurring as $\rho$ varies from $0.9$ to $0.72$, by assimilating noisy observations of $x_1$ and $x_2$ over the time window $[0,T=300]$. This is accomplished even though the EKF is equipped with a SINDy model that has been trained offline on trajectories corresponding to a single fixed parameter $\rho=0.92$. Therefore, the predictive model is initialised on trajectories leading to a fixed point, being unaware of the parametric dependency and the corresponding bifurcating behavior. The time step for integrating Eq.~\ref{eq:selkov} has been set to $\Delta t= 0.1$, which also serves as the sampling rate for the observations. Noisy observations have been generated by adding white noise to the signals, resulting in time series with a final SNR of $25$.

SINDy has been initially trained using the SQTLSQ algorithm with $\delta_r=0.05$ and $L=5\cdot 10^{-2}$, and a cubic polynomial functional library $\boldsymbol{\Theta}$. Results of the SINDy modelling are presented in Tab. \ref{tab:SELKOVcoefficients} in \ref{sec:appendixA}. In contrast to what is presented in Sec. \ref{sec:LV}, the function $\mathbf{f}$ describing the evolution of system dynamics, as identified offline by SINDy, differs from Eq.~\eqref{eq:selkov} in both the values of the parameters that weight the functional terms and in the number of functional terms included. In particular, SINDy introduces a quadratic contribution $x_1 x_2$ to model $\dot{x}_1=f_1(x_1,x_2)$, weighting it with the parameter $\iota$. It must be noted that increasing the SQTLSQ threshold to $L=10^{-1}$ would prevent SINDy from the identification of this additional term. However, modifying $L$ in this way would require a priori knowledge of the system dynamics, which is not always granted available. 
Assuming that this is one of those cases, we have investigated how the EKF performs when equipped with such, not correctly sparse, model. Here we have adopted an augmented state $\boldsymbol{\varkappa} = [\mathbf{x}^{\top},\tilde{\boldsymbol{\xi}}]^{\top}$, with $\mathbf{x}=[x_1,x_2]^{\top}$ and $\tilde{\boldsymbol{\xi}}=[\rho,\varsigma_1,\varsigma_2,\vartheta_1,\vartheta_2,\psi,\iota]$.

\begin{figure}[t!]
    \centering
    \subfloat[]{%
\includegraphics[width=145mm]{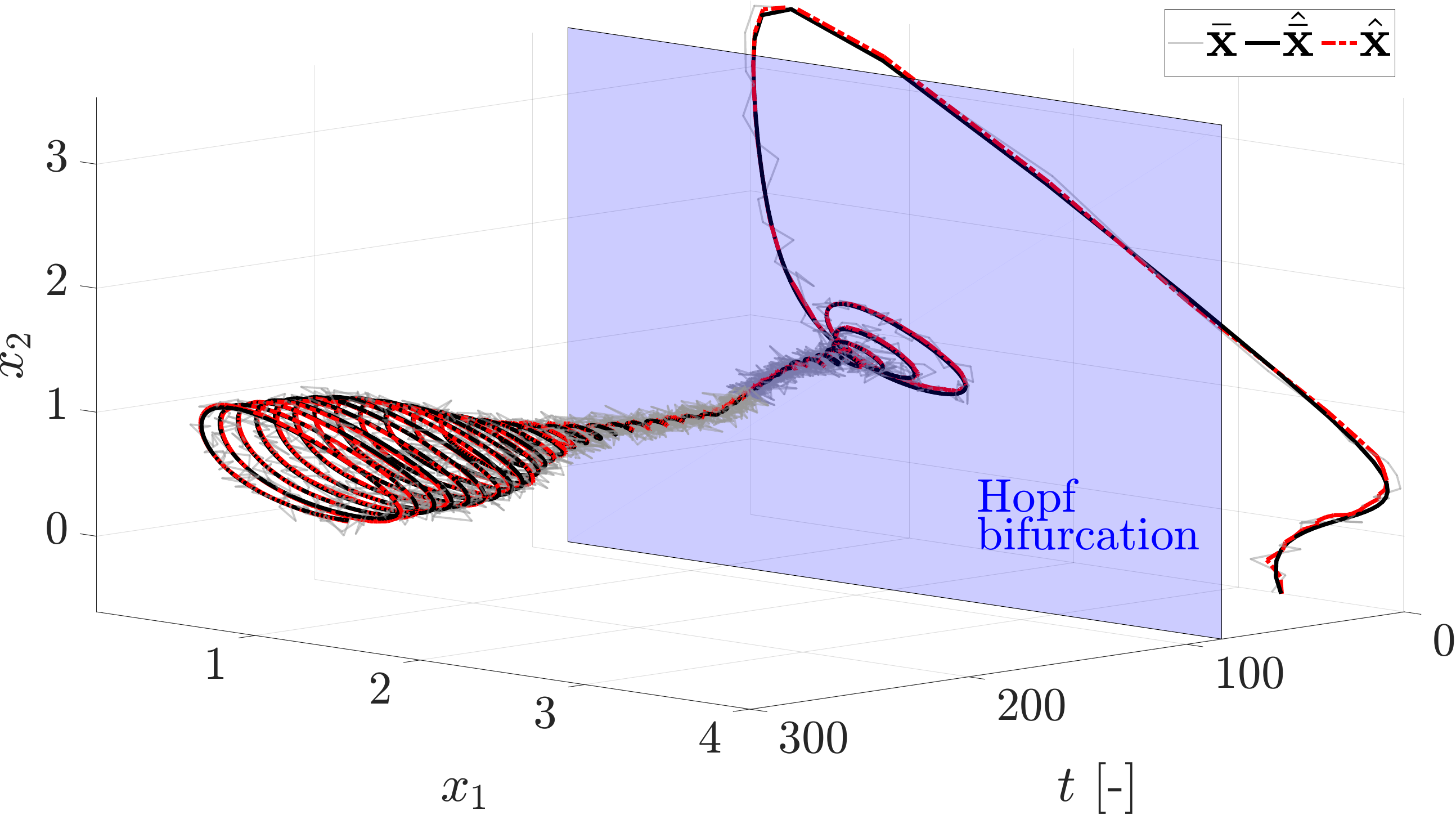}
\label{fig:Selkov_3D}}
\\
\subfloat[]{%
\includegraphics[width=145mm]{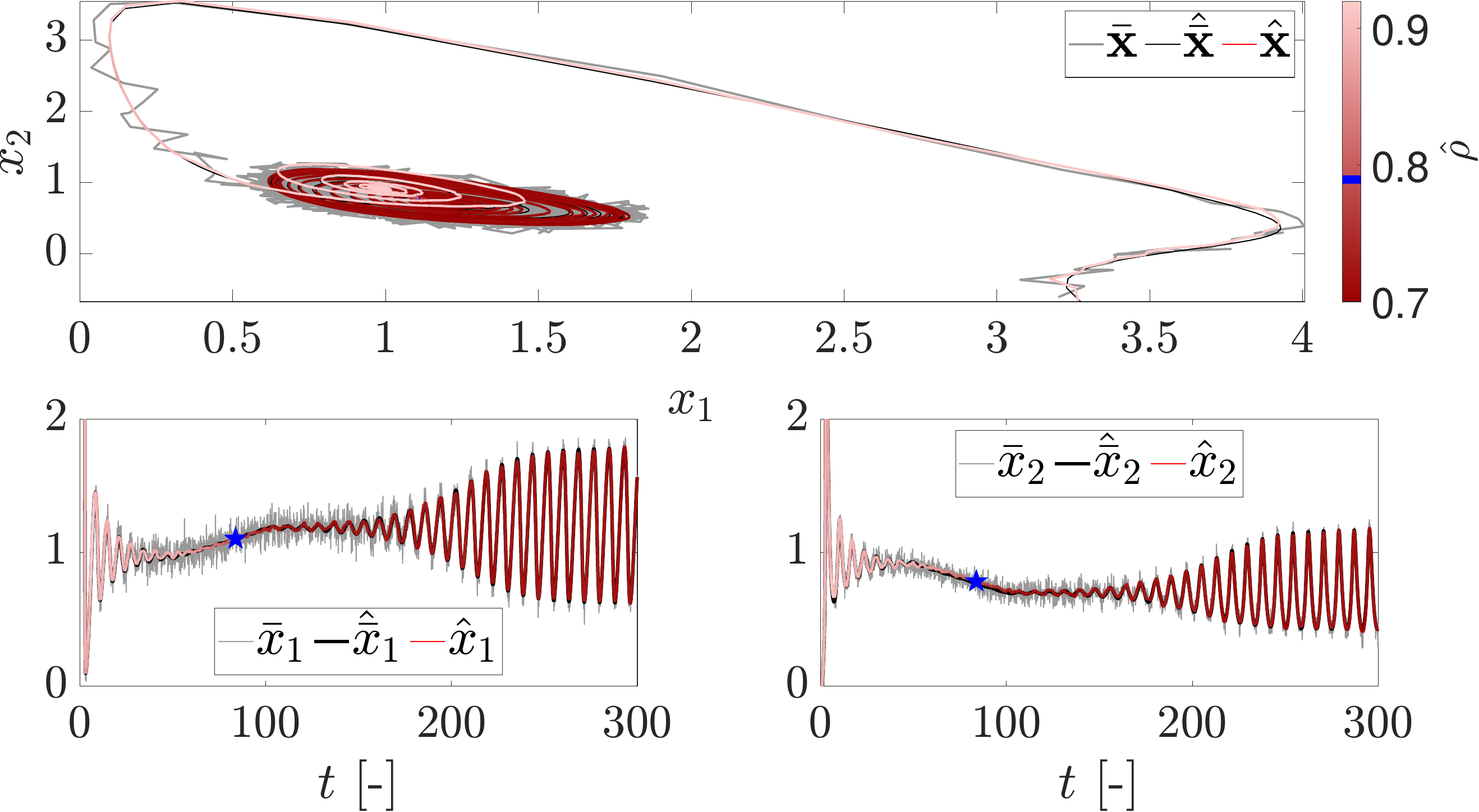}
\label{fig:Selkov_x12}}
\caption{Selkov model. In the top figure (a), the estimated states are plotted (red dotted line) against the assimilated noisy signals (gray solid lines) using a 3D state-space-time representation. Additionally, noise-free signals are reported in black solid lines. In the bottom plot (b), the projection of this 3D representation onto the 2D planes is reported. Colouring depending on the estimate $\hat{\rho}$ of the parameter ruling the two types of dynamics of the system (stable focus and stable limit cycles). The blue colour (further marked with a star in the $(t,x_1)$ and $(t,x_2)$ planes) indicates when the system undergoes a Hopf bifurcation.}
\label{fig:Selkov_x}
\end{figure}

\begin{figure}[t!]
    \centering
    \subfloat[$f_1$ parameter identification.]{%
\includegraphics[width=145mm]{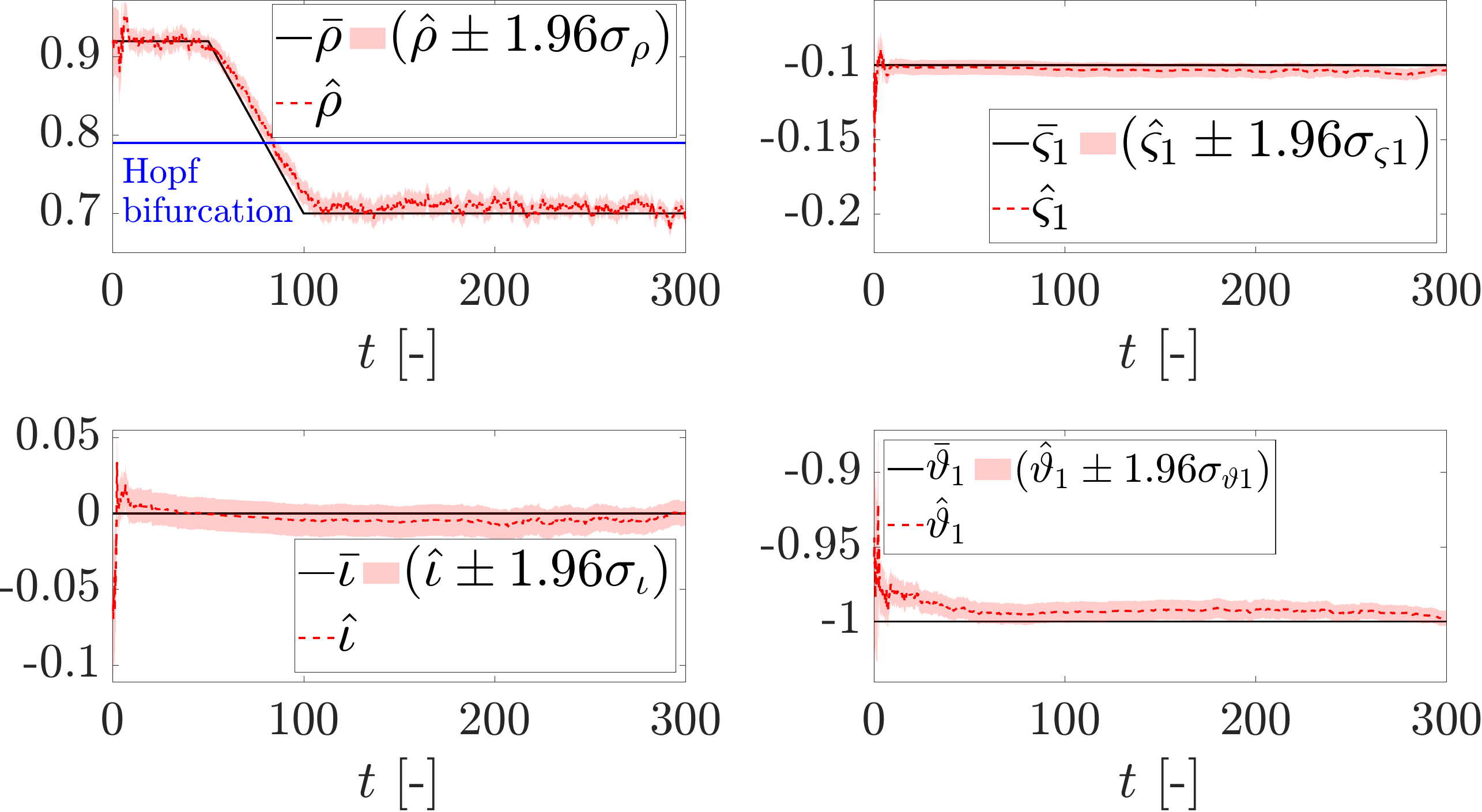}
\label{fig:Selkov_parameters_1}}
\\
\subfloat[$f_2$ parameter identification.]{%
\includegraphics[width=145mm]{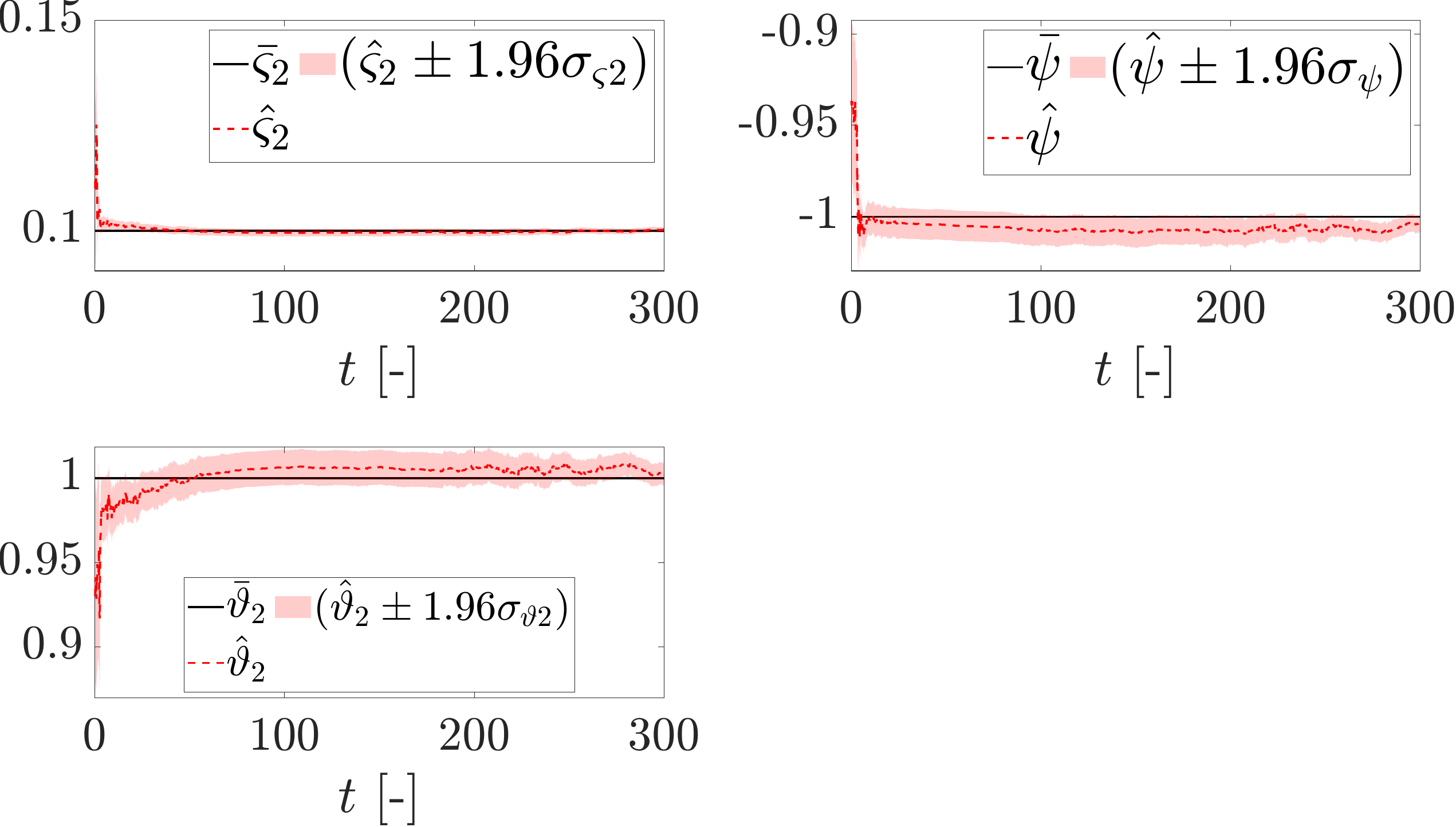}
\label{fig:Selkov_parameters_2}
}
\caption{Selkov model. The time-evolution of the estimated SINDy-model parameter are plotted (red dotted line) against the target time evolution (black solid line). The top figure reports the evolution of the parameters entering $\dot{x}_1=f_1(x_1,x_2)$; the bottom figure the evolution of the parameters entering $\dot{x}_2=f_2(x_1,x_2)$. The red shaded area represents $95\%$ confidence interval of the estimates, determined using the posterior covariance.}
\label{fig:Selkov_parameters}
\end{figure}

The filtering outcomes are presented in Figs. \ref{fig:Selkov_x} and \ref{fig:Selkov_parameters}, which respectively illustrate the identification results for the variables tracking the system dynamics, $x_1$ and $x_2$, and for the model parameters. In the second figure, it can be observed that the EKF is able not only to track the evolution of $\rho$ but also to improve the estimate of the other parameters. Specifically, the filter sets $\iota=0$, thereby restoring the correct sparsity of the system.

Notably, the recovery of the target parameter values has been achieved while tracking the system through two distinct dynamics regimes divided by a Hopf bifurcation. This is shown in Fig. \ref{fig:Selkov_3D}, which illustrates the evolution of $x_1$ and $x_2$ in $[0,T]$. Before the Hopf bifurcation at $t=83.8$, the system approaches a stable fixed point, although the location of this point shifts due to variations in the parameter $\rho$. After the bifurcation, the system dynamics change, resulting in orbits forming a stable limit cycle. Interestingly, the trajectories used to train the SINDy model do not include limit cycle orbits. However, due to the coupling with the filter, the resulting data-driven model is able to generalise to a completely different dynamics. Such an outcome may be difficult to achieve with other ML techniques, highlighting the significant value of the proposed approach. In Fig. \ref{fig:Selkov_x12}, the previously discussed 3D representation is projected onto the $(x_1,x_2)$, $(t,x_1)$, and $(t,x_2)$ planes for clarity. 

\subsection{MEMS arch}
\label{sec:MEMSarch}

\begin{figure}[t!]
\centering
\includegraphics[width=145mm]{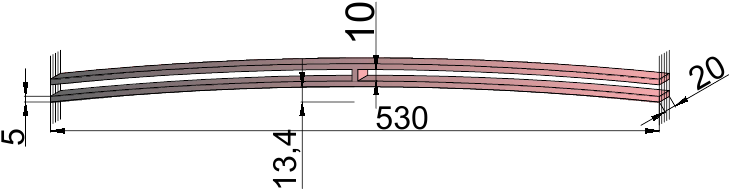}
{\caption{\footnotesize 
MEMS arch geometry. Dimensions in $\mu \text{m}$.}\label{fig:MEMS_arch}}
\end{figure}

The proposed procedure is finally applied to a MEMS clamped shallow arch, whose geometry is shown in Fig.~\ref{fig:MEMS_arch}. The device is made by polycrystalline silicon with density equal to $D=2330$ kg/m$^3$, Young modulus $E=1.670\cdot10^5$ MPa and Poisson ratio $\nu=0.22$. A linear elastic behaviour has been assumed to model its constitutive behaviour.

This structure has been extensively studied due to its intriguing dynamic response, which is characterised by a 1:2 internal resonance, and the emergence of a quasi-periodic regime as a result of of Neimark-Sacker (NS) bifurcations \cite{art:Giorgio23}. By adopting the two eigenmodes involved in the 1:2 resonance as master modes, a 2 Degrees of Freedom (DOFs) Reduced Order Model (ROM) was obtained in \cite{art:Giorgio21} via the static condensation approach \cite{art:Attilio19}. In this method, a stress manifold is generated by applying to the arch static loadings proportional to the inertia of the master modes. This manifold was thus fitted with a third order polynomial to get an explicit expression for the nonlinear restoring forces. A representation of the eigenmodes involved in the arch resonance can be found in \cite{art:Giorgio21}. The two resonance frequencies are $\omega_1=0.43416$ MHz and $\omega_2=0.86367$ MHz, respectively, while the $\mathcal{Q}$ factors are $\mathcal{Q}_1=500$ and $\mathcal{Q}_2=1000$. The connection between the $\mathcal{Q}$ factors and the damping coefficients will be specified in the following. The eigenmode amplitudes are represented by the modal coordinates $u_1$ and $u_2$, measured in $\mu \text{m}$, while their time derivative are denoted by $v_1$ and $v_2$. The following mass normalised ROM $\mathbf{f}$, here playing the role of physical model, expressed with $\tau = \omega_1 t$, is considered for the arch dynamics:
\begin{subequations}
    \begin{align}
        \dot{u}_1 = f_{1v}(v_1) = v_1, \\
        \dot{u}_2 = f_{2v}(v_2) = v_2, \\
        \dot{v}_1 = f_{1}(u_1,u_2,v_1) = -k_1 u_1 -\mu_1 v_1 -a_{12}u_1u_2 -b_{111}u_1^3 + B\text{cos}(\Omega \tau), \\
        \dot{v}_2 = f_{2}(u_1,u_2,v_2) = -k_2 u_2 -\mu_2 v_2 -a_{11}u_1^2 -b_{222}u_2^3,
    \end{align}
    \label{eq:ROM_MEMS_reduced}
\end{subequations}
\noindent where: $\mu_1 = 1/\mathcal{Q}_1$ and $\mu_2=\omega_2/(\mathcal{Q}_2\omega_1)$ are the damping coefficients of the two modes; $k_1$ and $k_2$ are the linear stiffness coefficients; $a_{12},a_{11}$ are the parameters ruling the quadratic contributions to the elastic force; $b_{111}$ and $b_{222}$ are parameters ruling the cubic contributions; $ B\text{cos}(\Omega \tau)$ is the harmonic forcing amplified by the load amplifier $B=0.0201 \mu\text{N}\mu\text{s}^2/\text{ng}$, that is applied by driving an alternated current through capacitor plates positioned near the arch. 

\begin{figure}[ht]
\centering
\includegraphics[width=145mm]{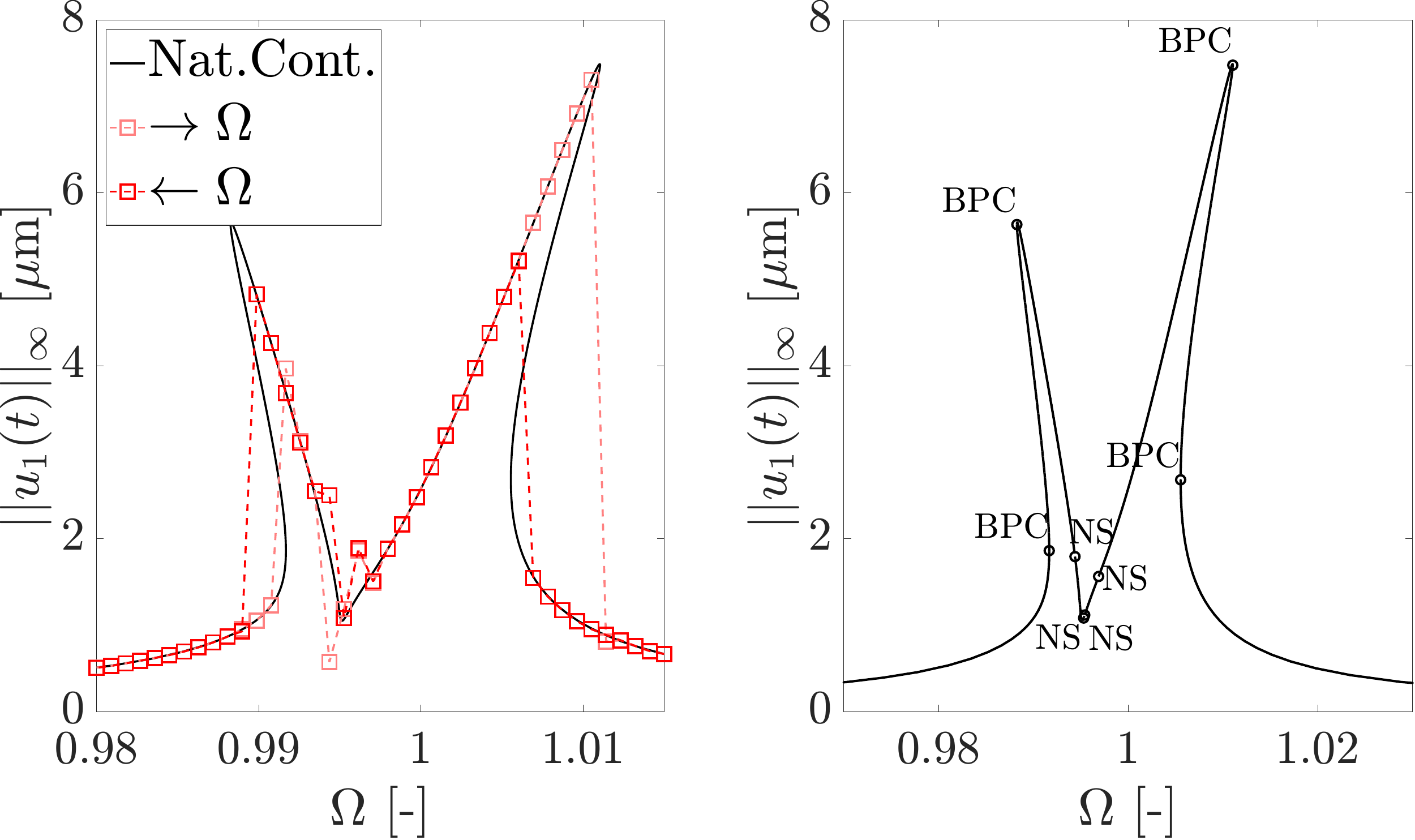}
{\caption{\footnotesize 
MEMS arch. On the left, comparison between the FRC obtained through numerical continuation (Nat. Cont.), and the FRCs obtained through time integration, both increasing ($\rightarrow \Omega$) and decreasing ($\rightarrow \Omega$) the excitation frequency $\Omega$. On the right, the FRC of the system is reported highlighting the Branch Point Cycles (BPC) and Neimark-Sacker (NS) bifurcations.}\label{fig:frf_tm}}
\end{figure}

Compared to the original formulation adopted in \cite{art:Giorgio21}, Eq.\eqref{eq:ROM_MEMS_reduced} retains the terms that govern the most significant aspects of the arch dynamics. Specifically, the internal resonance between the two modes is driven by the dependence of $\dot{v}_2$ on the quadratic term $a_{11}u_1^2$; consequently, we have also included $a_{12}u_1u_2$ to model $\dot{v}_1$, as this term accordingly derives from the elastic potential whenever $a_{11}u_1^2$ is present in the formulation of $\dot{v}_2$ \cite{art:Attilio19}. Cubic terms are also retained to capture the softening and hardening behaviour of the two peaks of the system Frequency Response Curve (FRC) shown in Fig. \ref{fig:frf_tm}. The other terms of \cite{art:Giorgio21} are not retained because they result only on a marginal shift of the bifurcation points. The FRC plots the maximum displacement $\tilde{u}_1$ of the first eigenmode at steady state against the non-dimensional excitation frequency $\Omega$. The adopted formulation enables the onset of the 1:2 internal resonance, NS bifurcations, and the presence of frequency combs between these bifurcation points. The figure illustrates both the precise locations of the bifurcation points (on the right), and a comparison between the FRC obtained through Natural Continuation methods with that obtained from the time integration of Eq.~\eqref{eq:ROM_MEMS_reduced} (on the left). The latter plot highlights the presence of frequency combs within the FRC region bounded by the NS bifurcation points, which are not observable when using a natural continuation method and assuming a single frequency periodic response. Natural continuation has been performed using \texttt{Matcont} \cite{art:MATCONT}.

First, we have recovered Eq.~\ref{eq:ROM_MEMS_reduced} via SINDy by observing the displacements and velocities of the arch eigenmodes, obtained by projecting the arch full field displacements and velocities. To this end, we have employed the SQTLSQ algorithm with a  regularisation strength of $\delta_r=0.05$ and threshold parameter $L=10^{-5}$. For the composition of the function library $\boldsymbol{\Theta}$, polynomial terms up to the third order have been used in accordance with Eq.~\eqref{eq:ROM_MEMS_reduced}, together with a cosine term describing the forcing condition. The composition of $\boldsymbol{\Theta}$ could have been inferred by developing the weak form of the dynamic equilibrium equation of the arch under large displacements and small strains \cite{art:Attilio19}. Quadratic terms naturally arise from the asymmetric shape of the arch, while cubic terms emerge from the clamped-clamped condition \cite{art:Attilio19,art:Giorgio23}. This type of reasoning, driven by the physics of the problem, may not be easily followed in other data-driven approaches to modelling dynamic systems such as recurrent neural networks, thus highlighting a significant advantage of the proposed procedure. The results of the identification are presented in the \ref{sec:appendixB} by specifying the coefficients ruling the linear, quadratic and cubic terms. This shows that SINDy has successfully recovered the dynamics of the system although no parametric variability has been shown to SINDy. Indeed, the matrix $\mathbf{X}$ collect trajectories for different forcing frequency $\Omega$, while the linear, quadratic, and cubic terms have been fixed to the values considered in \cite{art:Giorgio21}. 

After applying SINDy, we have used the EKF to adapt the identified dynamics in two scenarios, each posing distinct challenges. Both scenarios have been treated numerically, with no noise sources taken into account. The first scenario (Sec. \ref{sec:MEMS_sec_1}) involves updating the dynamic model $\mathbf{f}$, which describes the evolution of the two eigenmodes. Initially, this model is assumed to be identified via SINDy for a specific device, but it is then adapted to a MEMS with the same specifications however showing a different dynamic response. Since it is generally well assessed how to estimate the linear stiffness and damping terms, the focus has been on updating the estimate of the coefficients weighting the quadratic and cubic terms. Such an adaptation is commonly required in industrial settings to account for fabrication uncertainties, such as variations in the amount of overetch \cite{art:Mirzazadeh18}. Another important application is sensor self calibration \cite{art:Brossard20,proc:Saccani24}. Throughout its lifetime, a sensor may experience mechanical and thermal stress, that can cause unwanted shifts in its readings. Automatically compensating for errors is crucial to ensure the reliability of measurement systems, particularly in critical applications like structural health monitoring \cite{art:Giordano23}. 

The second scenario (Sec. \ref{sec:MEMS_sec_2}) involves recovering the coupling quadratic terms responsible for the internal resonance from data. This demonstrates how the procedure can enrich the dynamic description of the system while acquiring data. The only requirement is to include the relevant coefficients in the state vector, although initially set to an almost vanishing value. We believe that this scenario illustrates how the proposed approach could enable online identification of emerging contributions to the dynamics of a mechanical system.

\subsubsection{Estimation of the quadratic and cubic nonlinearities}
\label{sec:MEMS_sec_1}

We consider adapting the values of the quadratic and cubic coefficients identified by SINDy (see \ref{sec:appendixB}) to a different device which shares the same design but exhibits different dynamic behaviour which, for instance, is very typical of microstructures that are produced on identical wafers, but are characterized by slightly different geometrical features due to technological uncertainties. The adaptation is performed using the EKF, based on measurements of the eigenmode displacements and velocities. Specifically, this second device features quadratic and cubic coefficients that are $25\%$ higher than those of the original device. The state vector $\boldsymbol{\varkappa}=[\mathbf{x}^{\top},\boldsymbol{\zeta}^{\top}_{\Phi}]^{\top}$ includes both the dynamic variables $\mathbf{x}(t)=[u_1,u_2,v_1,v_2]^{\top}$ and the parameters $[a_{11},a_{12},b_{111},b_{222}]^{\top}$, which are treated as time-dependant variables.

The EKF equipped with the modified SINDy model has been applied over a time window of $[\tau_0,\tau_T]$, with $\tau_0=0$ and $\tau_T=5000$. The identification has been particularly challenging due to the small magnitude of the quadratic and cubic nonlinear coefficients. The system has been forced at $\Omega=0.9899$, corresponding to the first peak of the FRC. This forcing frequency has been chosen to enhance sensitivity to the cubic contribution to the internal forces.

\begin{figure}[ht!]
    \centering
    \subfloat[FRC modification during the identification (on top); FRC with uncertainty bounds at $t=T$ (on bottom).]{%
\includegraphics[width=145mm]{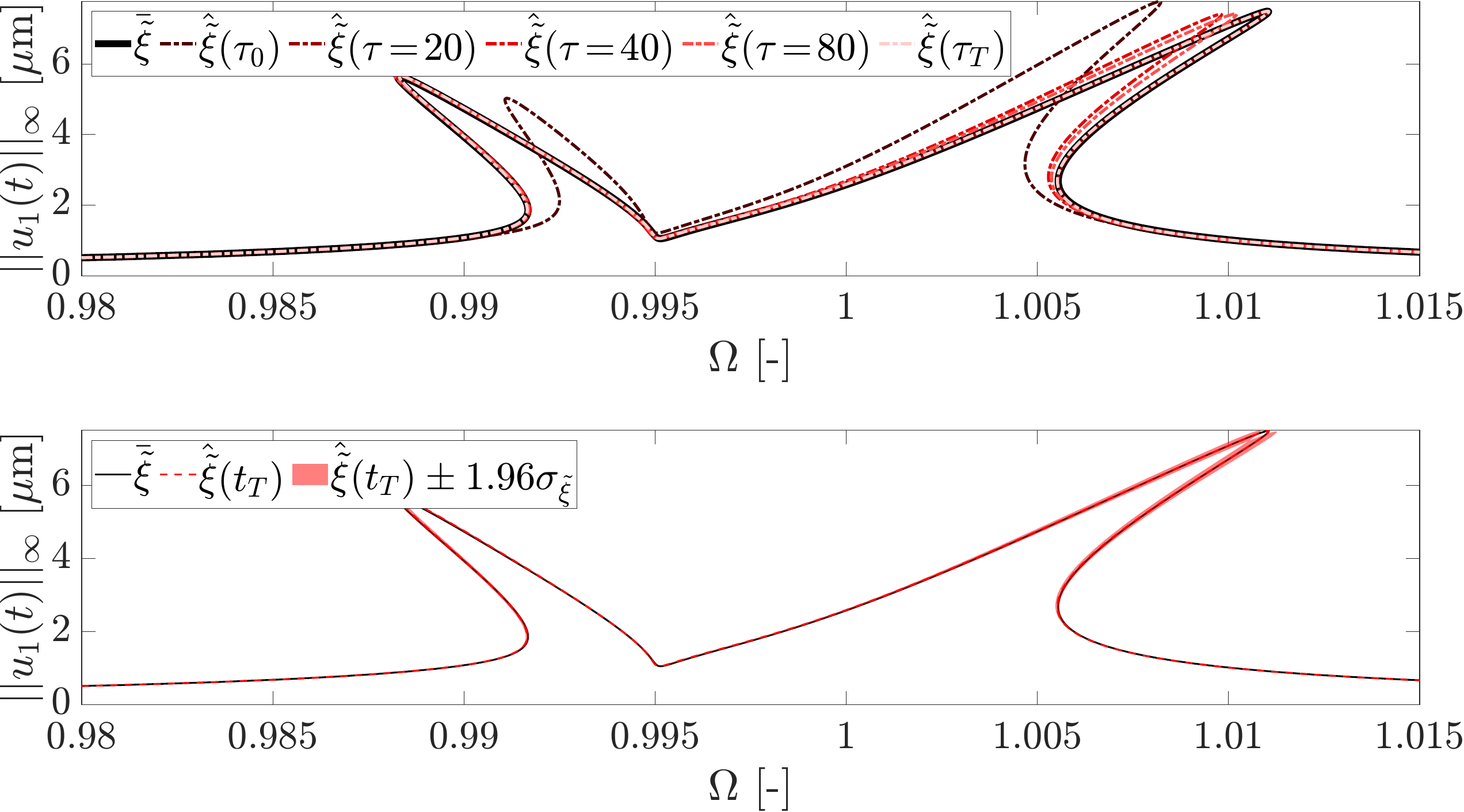}
\label{fig:MEMS_ev_frf_1}
} \\
    \subfloat[Parameter identification.]{%
\includegraphics[width=145mm]{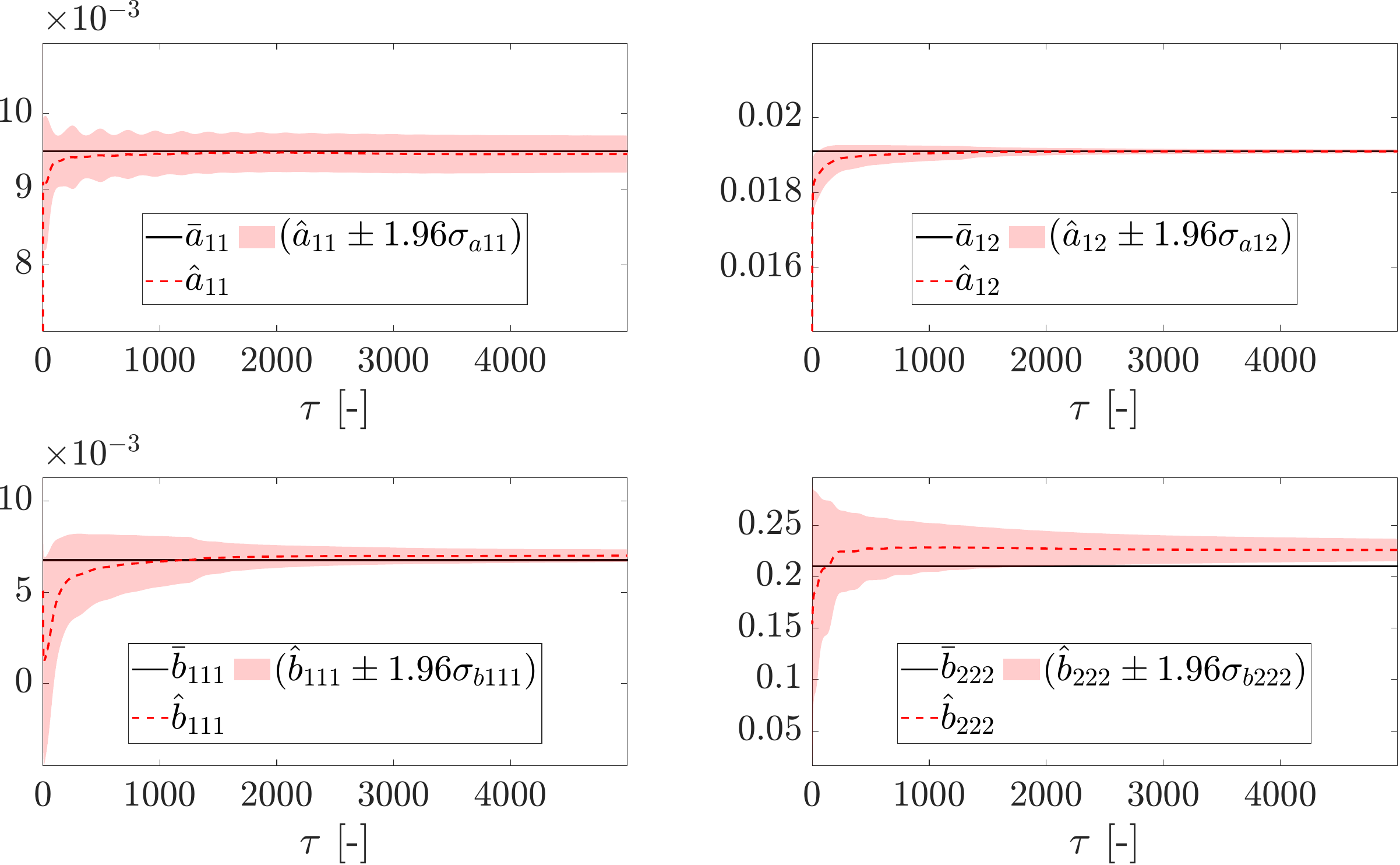}
\label{fig:MEMS_parameters1}}
\caption{MEMS arch, estimation of the quadratic and cubic nonlinearities. In the top figure, the FRCs corresponding to the models estimated during the identification (red palette dotted lines) are compared with the target FRC (black solid line). In the central figure,  The FRC identified at the end of the analysis (red dotted line) is plot against the target FRC (black solid line), with the red shaded area representing the $95\%$ confidence intervals of the $\tilde{\boldsymbol{\xi}}$ estimates, determined using the posterior covariance. In the bottom figure, the time-evolution of the estimated SINDy-model parameter are plotted (red dotted line) against the target time evolution (black solid line).}
\label{fig:MEMS_case1}
\end{figure}

Results are shown in Fig. \ref{fig:MEMS_case1}, showing excellent performance. Notably, key information about the dynamic response of the device are contained in the FRCs reported in Fig. \ref{fig:MEMS_ev_frf_1}, where the adimensionalised excitation frequency $\Omega$ is plotted on the horizontal axis and $\tilde{u}_1$ on the vertical axis. Similar figures can be realised with $\tilde{u}_2$ on the vertical axis, leading to similar comments to those made for $\tilde{u}_1$.

On top of Fig. \ref{fig:MEMS_ev_frf_1}, it is shown how the FRC changes during the estimation. Specifically, it can be seen how far the FRC, initialised (at $t_0$) with the $\tilde{\boldsymbol{\xi}}$ values determined by SINDy on the reference device, is from the target ones, and how quickly the EKF updates improve the estimations of $\tilde{\boldsymbol{\xi}}$. A satisfactory outcome is achieved already at $\tau=20$. This rapid response highlights the potential for online application of the proposed procedure, making it suitable for handling evolving environmental conditions that may affect the device dynamic response, or for responding to unexpected changes in the arch mechanics, such as damage from impacts. At the bottom of Fig. \ref{fig:MEMS_ev_frf_1}, the estimated FRC is shown along with its corresponding uncertainty bounds. These bounds are generally very tight but become slightly larger close to the FRC peaks. This widening is more pronounced near the peak corresponding at $\Omega=1.012$, likely due to the chosen excitation frequency being related to the first peak.

Fig. \ref{fig:MEMS_parameters1} illustrates the time evolution of the parameter estimation for $\tilde{\boldsymbol{\xi}}$. It can be noted that the estimated $\hat{b}_{222}$ does not exactly match the target value of $b_{222}$. However, we argue that the filter is unlikely to improve further, as the system dynamics prediction, as shown by the FRCs, is already highly accurate. In contrast, the estimation of the other cubic term coefficient $\hat{b}_{111}$, as well as the of quadratic terms, is excellent, with a remarkable reduction in estimation uncertainties over the course of the analysis.

\subsubsection{Discovering the internal resonance}
\label{sec:MEMS_sec_2}

We now consider the case in which $a_{11}$ is set to small enough ($1\cdot10^{-9}$) such that it does not induce internal resonance in the arch. It is noted that $a_{11}$ can not be set exactly to zero, as this would eliminate sensitivity to the parameter when computing the Jacobian matrix.
We have thus included $a_{11}$ in the state vector $\boldsymbol{\varkappa}=[\mathbf{x},\tilde{\boldsymbol{\xi}}]$ among the coefficients $\tilde{\boldsymbol{\xi}}=[a_{11},a_{12}]^{\top}$. Unlike the previous discussion, we have assumed that $b_{111}$ and $b_{222}$ are known to simplify the system identification, while we have still considered to underestimate the value of $a_{12}$ by $25\%$.

\begin{figure}[t!]
\centering
\includegraphics[width=145mm]{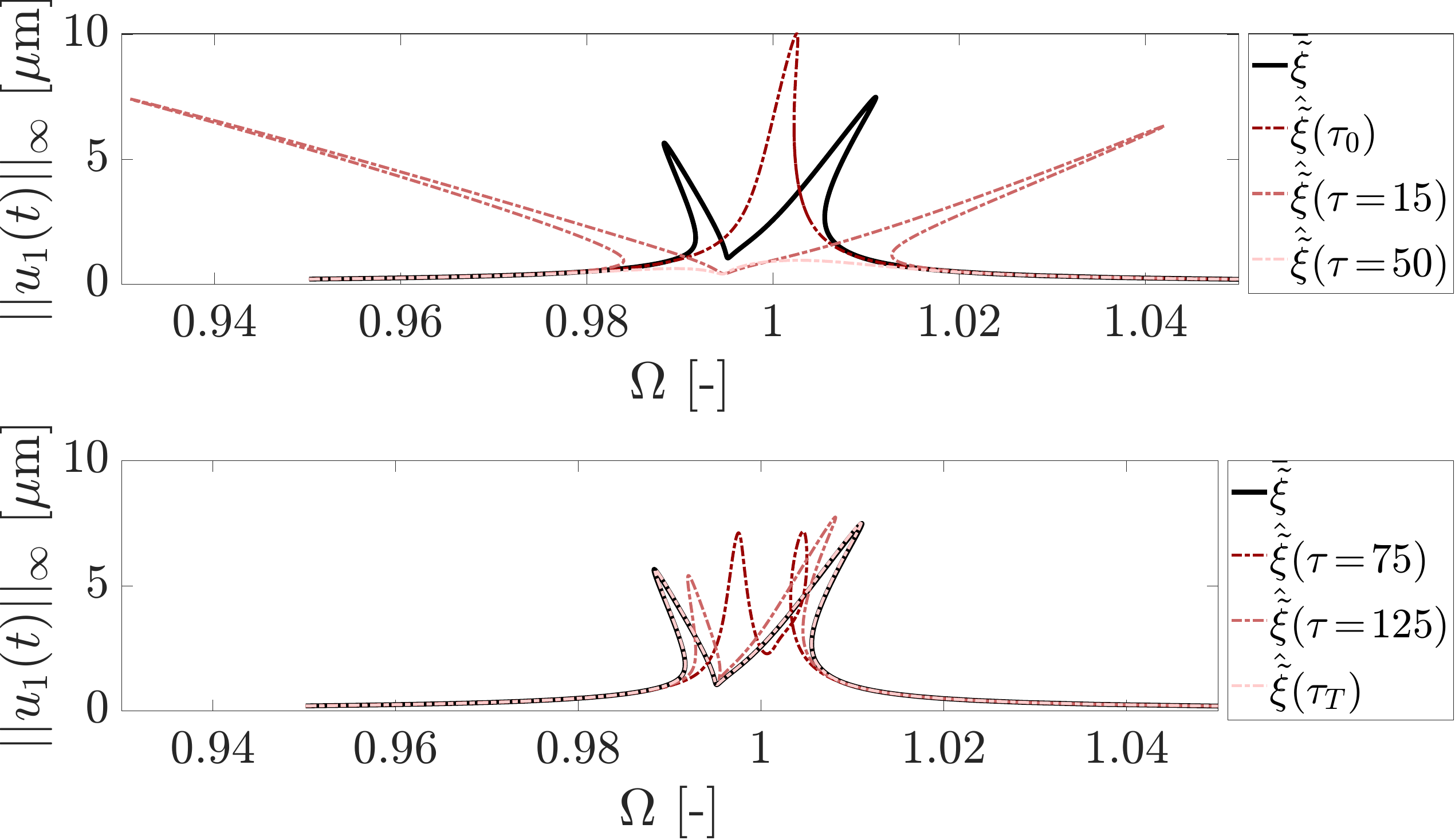}
{\caption{\footnotesize 
MEMS arch, discovering the internal resonance. The FRCs corresponding to the models estimated during the identification (red palette dotted lines) are compared with the target FRC (black solid line).}\label{fig:MEMS_case2}}
\end{figure}

Results are shown in Fig. \ref{fig:MEMS_case2}. At $t=0$, the FRC of the system, based on the initial estimates $\tilde{\boldsymbol{\xi}}(\tau=0)$, exhibits a single peak, as the value of $a_{11}$ is too small to induce the internal resonance. As the quadratic terms estimates are updated and the corresponding FRCs are plotted, the system dynamics gradually converge to the true behaviour. Specifically, at $\tau=15$, the initial peak levels out, and the appearance of two small amplified hills indicates that the identified system exhibits some interaction between the two eigenmodes. By $\tau=50$, the internal resonance becomes marked, with two peaks showing significant hardening and softening behaviours. At $\tau=125$, the identified system has the correct FRC shape, thus almost perfectly matching the target $\bar{\boldsymbol{\zeta}}_{\Phi}$ at the end of the analysis ($\tau_T$).

\section{Conclusions}
\label{sec:conclusions}
A combined approach exploiting the extended Kalman Filter (EKF) and the Sparse Identification of Nonlinear Dynamics (SINDy) technique has been adopted for online data assimilation and adaptation of nonlinear dynamics systems. SINDy is first used in a preliminary stage to identify a dynamic model, avoiding biases from incorrect assumptions. Second, the EKF is employed to data assimilation, directly updating the SINDy parameters. This procedure has shown the ability to track time-evolving systems, performing joint state-parameter estimation even for SINDy coefficients not accounted for during the relevant initial SINDy training. Notably, the EKF has been able to track systems that have undergone unseen bifurcations in their dynamics behaviour.

The promising results suggest that this procedure could be applied to a wide range of time-varying problems in dynamics. Future work could extend the method to partially unobserved system by incorporating time-delay embedding techniques \cite{art:Brunton17}, similarly to \cite{art:CMAME24}. Conversely, in cases where the number of observations exceeds the number of dominant dynamic modes, autoencoders architectures could be employed. To further enhance the robustness of this approach, automatic procedures for the EKF tuning will be also addressed in future research.

\section*{Code and data accessibility}
The source code of the proposed method is made available from the GitHub repository:\\
\texttt{https://github.com/ContiPaolo/EKF-SINDy} \cite{EKFSINDY_repo}. 

\section*{Acknowledgments}
The authors would like to thank Alessio Colombo (Politecnico di Milano) for the invaluable insights on numerical continuation methods. LR is supported by the Joint Research Platform ``Sensor sysTEms and Advanced Materials" (STEAM) between Politecnico di Milano and STMicroelectronics. AM 
acknowledges the project FAIR (Future Artificial Intelligence Research), funded by the NextGenerationEU program within the PNRR-PE-AI scheme (M4C2, Investment 1.3, Line on Artificial Intelligence). AM and PC acknowledge the PRIN 2022 Project “Numerical approximation of uncertainty quantification problems for PDEs by multi-fidelity methods (UQ-FLY)” (No. 202222PACR), funded by the European Union - NextGenerationEU.
AF acknowledges the PRIN 2022 Project “DIMIN- DIgital twins of nonlinear MIcrostructures with iNnovative model-order-reduction strategies” 
(No. 2022XATLT2) funded by the European Union - NextGenerationEU.

\appendix
\section{}
\label{sec:appendixA}

In the following, we report the comparison between the (reference) model parameters used for the training of SINDy, weighting the relevant functional terms, and the parameters identified by SINDy for the Lokta-Volterra model (Tab. \ref{tab:LVcoefficients}), the Selkov model (Tab. \ref{tab:SELKOVcoefficients}), and the MEMS arch (Tab. \ref{tab:ROMcoefficients}).

\paragraph{Lokta Volterra model}

\begin{center}
\begin{tabular}{cccc|cc}
\toprule
& & \multicolumn{2}{c|}{$\dot{x}_1 = f_1(x_1,x_2)$} &  \multicolumn{2}{c}{$\dot{x}_2 = f_2(x_1,x_2)$} \\
\midrule
\textbf{Coefficients} & Relevant term & \textbf{Reference} &  \textbf{SINDy} & \textbf{Reference} &  \textbf{SINDy} \\
$a$ & $x_1$ & $1.000$ & $1.000$ & $0$ & $0$\\
$b$ & $x_1x_2$ & $-0.1000$ & $-0.1000$ & $0$ & $0$\\
$c$ & $x_2$ & $0$ & $0$ & $-1.500$ & $-1.500$\\
$d$ & $x_1 x_2$ & $0$ & $0$ & $7.500\cdot10^{-2}$ & $7.500\cdot10^{-2}$\\
\bottomrule
\end{tabular}
\captionof{table}{Lokta Volterra model. The ``Reference" columns report the parameters used for the training of SINDy; in the ``SINDy" columns, we list the parameter values identified by SINDy\label{tab:LVcoefficients}.}
\end{center}

\paragraph{Selkov model}

\begin{center}
\begin{tabular}{cccc|cc}
\toprule
& & \multicolumn{2}{c|}{$\dot{x}_1 = f_1(x_1,x_2)$} &  \multicolumn{2}{c}{$\dot{x}_2 = f_2(x_1,x_2)$} \\
\midrule
\textbf{Coefficients} & Relevant term & \textbf{Reference} &  \textbf{SINDy} & \textbf{Reference} &  \textbf{SINDy} \\
$\rho$ & $1.0$ & $0.9200$ & $0.9234$ & $0$ & $0$\\
$\varsigma_1$ & $x_1$ & $-0.1000$ & $-9.389\cdot 10^{-2}$ & $0$ & $0$\\
$\varsigma_2$ & $x_1$ & $0$ & $0$ & $0.1000$ & $0.1082$\\
$\psi$ & $x_2$ & $0$ & $0$ & $-1.000$ & $-0.9343$ \\
$\iota$ & $x_1x_2$ & $0$ & $ -7.641\cdot 10^{-2}$ & $0$ & $0$\\
$\vartheta_1$ & $x_1 x^2_2$ & $-1.000$ & $-0.9294$ & $0$ & $0$\\
$\vartheta_2$ & $x_1 x^2_2$ & $0$ & $0$ & $1.000$ & $0.9185$\\
\bottomrule
\end{tabular}
\captionof{table}{Selkov model. The ``Reference" columns report the parameters used for the training of SINDy; in the ``SINDy" columns, we list the parameter values identified by SINDy\label{tab:SELKOVcoefficients}.}
\end{center}

\paragraph{MEMS arch}

\begin{center}
\begin{tabular}{cccc|cc}
\toprule
& & \multicolumn{2}{c|}{$\dot{v}_1 = f_1(u_1,u_2,v_1)$} &  \multicolumn{2}{c}{$\dot{v}_2 = f_2(u_1,u_2,v_2)$} \\
\midrule
\textbf{Coefficients} & Relevant term & \textbf{Reference} &  \textbf{SINDy} & \textbf{Reference} &  \textbf{SINDy} \\
\midrule
$k_1$ & $u_1$ & $1.000$ & $1.000$ & $0$ & $0$ \\
$k_2$ & $u_2$ & $0$ & $0$ & $3.957$ & $3.957$ \\
$\mu_1$ & $v_1$ & $2.000\cdot 10^{-3}$ & $2.000\cdot 10^{-3}$ & $0$ & $0$ \\
$\mu_2$ & $v_2$ & $0$ & $0$ & $2.000\cdot 10^{-3}$ & $1.999\cdot 10^{-3}$ \\
$a_{11}$ & $u_1^2$ & $0$ & $0$ & $7.125\cdot 10^{-3}$ & $7.124\cdot 10^{-3}$ \\
$a_{12}$ & $u_1u_2$ & $1.433\cdot 10^{-2}$ & $1.432\cdot 10^{-2}$ & $0$ & $0$ \\
$b_{111}$ & $u_1^3$ & $5.074\cdot 10^{-5}$ & $5.080\cdot 10^{-5}$ & $0$ & $0$ \\
$b_{222}$ & $u_2^3$ & $0$ & $0$ & $1.575\cdot 10^{-3}$ & $1.575\cdot 10^{-3}$\\
- & $B\text{cos}(\Omega\tau)$ & $1.000$ & $0.9999$ & $0$ & $0$ \\   
\bottomrule
\end{tabular}
\captionof{table}{MEMS arch. Coefficients weighting the linear, quadratic and cubic contribution to the elastic force. The ``Reference" column reports the reduced order model coefficients determined in \cite{art:Giorgio21} and used for the training of SINDy; in the ``SINDy" columns, we list the parameter values identified by SINDy\label{tab:ROMcoefficients}.}
\end{center}

\section{}
\label{sec:appendixB}

\scriptsize{
\begin{minipage}[h]{0.28\textwidth}
\label{tab:EKFtuning_1}
\begin{tabular}{lr}
\toprule
\multicolumn{2}{c}{\textbf{Lokta Volterra}} \\
\midrule
    $x_1, x_2$ & $10^{-3}$ \\
    $a$ & $10^{-4}$\\
    $b,c,d$ & $10^{-7}$\\
\midrule
\multicolumn{2}{c}{\textbf{Selkov}}\\
\midrule
    $x_1, x_2$ & $10^{-8}$\\
    $\rho, \varsigma_1,\vartheta_1,\vartheta_2,\psi,\iota$ & $10^{-3}$\\
    $\varsigma_2$ & $10^{-4}$\\
\midrule
\multicolumn{2}{c}{\textbf{MEMS arch}}\\
\midrule
\multicolumn{2}{c}{Quadratic and cubic}\\
\multicolumn{2}{c}{nonlinearities}\\
\midrule
    $u_1, u_2, v_1, v_2$ & $10^{-8}$\\
    $a_{11}$ & $5\cdot 10^{-4}$ \\
    $a_{12}$ & $10^{-4}$ \\
    $b_{111}$ & $10^{-5}$ \\
    $b_{222}$ & $5\cdot 10^{-3}$ \\
 \midrule
\multicolumn{2}{c}{Discovering internal}\\
\multicolumn{2}{c}{resonance}\\
    $u_1, u_2, v_1, v_2$ & $10^{-7}$\\
    $a_{11},a_{12}$ & $10^{-4}$ \\
\bottomrule
\end{tabular}
\captionof{table}{Initial covariance $\mathbf{P}(t_0)$, diagonal terms. In the left column, we specify the augmented state quantities to which these values are related.}
\end{minipage}%
$\quad \quad$
\begin{minipage}[h]{0.28\textwidth}
\label{tab:EKFtuning_2}
\begin{tabular}{lr}
\toprule
\multicolumn{2}{c}{\textbf{Lokta Volterra}} \\
\midrule
    $x_1, x_2$ & $10^{-3}$ \\
    $a$ & $5\cdot 10^{-5}$\\
    $b$ & $10^{-8}$\\
    $c$ & $10^{-14}$\\
    $d$ & $8\cdot 10^{-8}$\\
     \midrule
\multicolumn{2}{c}{\textbf{Selkov}}\\
\midrule
    $x_1, x_2$ & $10^{-6}$\\
    $\rho, \varsigma_1,\varsigma_2$ & $10^{-14}$\\
    $\vartheta_1, \vartheta_2,\psi,\iota$ & $10^{-12}$\\
    \midrule
\multicolumn{2}{c}{\textbf{MEMS arch}}\\
\midrule
\multicolumn{2}{c}{Quadratic and cubic}\\
\multicolumn{2}{c}{nonlinearities}\\
\midrule
    $u_1, u_2$ & $10^{-7}$\\
    $v_1$ & $10^{-5}$\\
    $v_2$ & $8 \cdot 10^{-6}$\\
    $a_{11}$ & $10^{-7}$ \\
    $a_{12}$ & $10^{-12}$ \\
    $b_{111}$ & $10^{-14}$ \\
    $b_{222}$ & $10^{-16}$ \\
\midrule
\multicolumn{2}{c}{Discovering internal}\\
\multicolumn{2}{c}{resonance}\\
    $u_1, u_2$ & $10^{-8}$\\
    $v_1, v_2$ & $10^{-5}$\\
    $a_{11}$ & $10^{-7}$ \\
    $a_{12}$ & $10^{-9}$ \\
\bottomrule
\end{tabular}
\captionof{table}{Process noise $\mathbf{Q}$, diagonal terms. In the left column, we specify the augmented state quantities to which these values are related.}
\end{minipage}%
$\quad \quad$
\noindent
\begin{minipage}[h]{0.28\textwidth}
\label{tab:EKFtuning_3}
\begin{tabular}{cr}
\toprule
\multicolumn{2}{c}{\textbf{Lokta Volterra}} \\
\midrule
    $x_1, x_2$ & $1.0$ \\
\midrule
\multicolumn{2}{c}{\textbf{Selkov}}\\
\midrule
    $x_1, x_2$ & $5\cdot10^{-4}$\\
\midrule
\multicolumn{2}{c}{\textbf{MEMS arch}}\\
\midrule
\multicolumn{2}{c}{Quadratic and cubic}\\
\multicolumn{2}{c}{nonlinearities}\\
\midrule
    $u_1, u_2$ & $10^{-1}$\\
    $v_1, v_2$ & $10^{-3}$\\
\midrule
\multicolumn{2}{c}{Discovering internal}\\
\multicolumn{2}{c}{resonance}\\
    $u_1, u_2$ & $10^{-1}$\\
    $v_1, v_2$ & $10^{-3}$\\
\bottomrule
\end{tabular}
\captionof{table}{Measurement noise $\mathbf{R}$, diagonal terms. In the left column, we specify the observed quantities to which these values are related.}
\end{minipage}%
}


\scriptsize \bigskip

\bibliographystyle{elsarticle-num}
\bibliography{references}

\begin{thebibliography}{10}
\expandafter\ifx\csname url\endcsname\relax
  \def\url#1{\texttt{#1}}\fi
\expandafter\ifx\csname urlprefix\endcsname\relax\def\urlprefix{URL }\fi
\expandafter\ifx\csname href\endcsname\relax
  \def\href#1#2{#2} \def\path#1{#1}\fi

\bibitem{book:Farrar13}
C.~Farrar, K.~Worden, Structural Health Monitoring: A Machine Learning Perspective, Wiley, Hoboken, NJ, 2013.
\newblock \href {https://doi.org/10.1002/9781118443118} {\path{doi:10.1002/9781118443118}}.

\bibitem{art:Fonzi20}
N.~Fonzi, S.~L. Brunton, U.~Fasel, Data-driven nonlinear aeroelastic models of morphing wings for control, Proceedings of the Royal Society A: Mathematical, Physical and Engineering Sciences 476~(2239) (2020) 20200079.
\newblock \href {https://doi.org/10.1098/rspa.2020.0079} {\path{doi:10.1098/rspa.2020.0079}}.

\bibitem{book:Corigliano18}
A.~Corigliano, R.~Ardito, C.~Comi, A.~Frangi, A.~Ghisi, S.~Mariani, Accelerometers, John Wiley \& Sons, Ltd, 2018, Ch.~4, pp. 91--108.
\newblock \href {https://doi.org/10.1002/9781119053828.ch4} {\path{doi:10.1002/9781119053828.ch4}}.

\bibitem{art:Wagg20}
D.~J. Wagg, K.~Worden, R.~J. Barthorpe, P.~Gardner, {Digital Twins: State-of-the-Art and Future Directions for Modeling and Simulation in Engineering Dynamics Applications}, ASCE-ASME Journal of Risk and Uncertainty in Engineering Systems, Part B: Mechanical Engineering 6~(3) (2020) 030901.
\newblock \href {https://doi.org/10.1115/1.4046739} {\path{doi:10.1115/1.4046739}}.

\bibitem{art:Matt24}
M.~Torzoni, M.~Tezzele, S.~Mariani, A.~Manzoni, K.~E. Willcox, A digital twin framework for civil engineering structures, Computer Methods in Applied Mechanics and Engineering 418 (2024) 116584.
\newblock \href {https://doi.org/10.1016/j.cma.2023.116584} {\path{doi:10.1016/j.cma.2023.116584}}.

\bibitem{art:Kapteyn21}
M.~G. Kapteyn, J.~V.~R. Pretorius, K.~E. Willcox, A probabilistic graphical model foundation for enabling predictive digital twins at scale, Nature Computational Science 1 (2021) 337--347.
\newblock \href {https://doi.org/10.1038/s43588-021-00069-0} {\path{doi:10.1038/s43588-021-00069-0}}.

\bibitem{art:Kalman60}
R.~E. Kalman, A new approach to linear filtering and prediction problems, Journal of Basic Engineering 82~(1) (1960) 35--45.
\newblock \href {https://doi.org/10.1115/1.3662552} {\path{doi:10.1115/1.3662552}}.

\bibitem{book:Simon06_5}
D.~Simon, The discrete--time {Kalman} filter, John Wiley \& Sons, Ltd, 2006, Ch.~5, pp. 121--148.
\newblock \href {https://doi.org/10.1002/0470045345.ch5} {\path{doi:10.1002/0470045345.ch5}}.

\bibitem{proc:Wan00}
E.~Wan, R.~Van Der~Merwe, The unscented {Kalman} filter for nonlinear estimation, in: Proceedings of the IEEE 2000 Adaptive Systems for Signal Processing, Communications, and Control Symposium (Cat. No.00EX373), 2000, pp. 153--158.
\newblock \href {https://doi.org/10.1109/ASSPCC.2000.882463} {\path{doi:10.1109/ASSPCC.2000.882463}}.

\bibitem{art:Mariani07}
S.~Mariani, A.~Ghisi, Unscented {Kalman} filtering for nonlinear structural dynamics, Nonlinear Dynamics 49 (2007) 131--150.
\newblock \href {https://doi.org/10.1007/s11071-006-9118-9} {\path{doi:10.1007/s11071-006-9118-9}}.

\bibitem{art:RUENG24}
L.~Rosafalco, S.~E. Azam, S.~Mariani, A.~Corigliano, System identification via unscented {Kalman} filtering and model class selection, ASCE-ASME Journal of Risk and Uncertainty in Engineering Systems, Part A: Civil Engineering 10~(1) (2024) 04023063.
\newblock \href {https://doi.org/10.1061/AJRUA6.RUENG-1085} {\path{doi:10.1061/AJRUA6.RUENG-1085}}.

\bibitem{art:Novoa24}
A.~N\'{o}voa, A.~Racca, L.~Magri, Inferring unknown unknowns: Regularized bias-aware ensemble {Kalman} filter, Computer Methods in Applied Mechanics and Engineering 418 (2024) 116502.
\newblock \href {https://doi.org/10.1016/j.cma.2023.116502} {\path{doi:10.1016/j.cma.2023.116502}}.

\bibitem{art:Song20}
M.~Song, R.~Astroza, H.~Ebrahimian, B.~Moaveni, C.~Papadimitriou, Adaptive kalman filters for nonlinear finite element model updating, Mechanical Systems and Signal Processing 143 (2020) 106837.
\newblock \href {https://doi.org/10.1016/j.ymssp.2020.106837} {\path{doi:10.1016/j.ymssp.2020.106837}}.

\bibitem{art:Gres25}
S.~Gre{\'s}, M.~D{\"o}hler, V.~K. Dertimanis, E.~N. Chatzi, Subspace-based noise covariance estimation for kalman filter in virtual sensing applications, Mechanical Systems and Signal Processing 222 (2025) 111772.
\newblock \href {https://doi.org/10.1016/j.ymssp.2024.111772} {\path{doi:10.1016/j.ymssp.2024.111772}}.

\bibitem{art:Bilgin25}
N.~Bilgin, A.~Olivier, Joint bayesian estimation of process and measurement noise statistics in nonlinear kalman filtering, Mechanical Systems and Signal Processing 223 (2025) 111836.
\newblock \href {https://doi.org/10.1016/j.ymssp.2024.111836} {\path{doi:10.1016/j.ymssp.2024.111836}}.

\bibitem{art:Cuomo24}
S.~Cuomo, M.~{De Rosa}, F.~Piccialli, L.~Pompameo, Railway safety through predictive vertical displacement analysis using the pinn-ekf synergy, Mathematics and Computers in Simulation 223 (2024) 368--379.
\newblock \href {https://doi.org/10.1016/j.matcom.2024.04.026} {\path{doi:10.1016/j.matcom.2024.04.026}}.

\bibitem{art:Impraimakis24}
M.~Impraimakis, Deep recurrent-convolutional neural network learning and physics kalman filtering comparison in dynamic load identification, Structural Health Monitoring 0~(0) (0) 14759217241262972.
\newblock \href {https://doi.org/10.1177/14759217241262972} {\path{doi:10.1177/14759217241262972}}.

\bibitem{art:SDEE23}
F.~Gatti, L.~Rosafalco, G.~Colombera, S.~Mariani, A.~Corigliano, Multi-storey shear type buildings under earthquake loading: Adversarial learning-based prediction of the transient dynamics and damage classification, Soil Dynamics and Earthquake Engineering 173 (2023) 108141.
\newblock \href {https://doi.org/10.1016/j.soildyn.2023.108141} {\path{doi:10.1016/j.soildyn.2023.108141}}.

\bibitem{art:Joseph24}
H.~Joseph, G.~Quaranta, B.~Carboni, W.~Lacarbonara, Deep learning architectures for data-driven damage detection in nonlinear dynamic systems under random vibrations, Nonlinear Dynamics 223 (2024).
\newblock \href {https://doi.org/10.1007/s11071-024-10270-1} {\path{doi:10.1007/s11071-024-10270-1}}.

\bibitem{proc:Coskun17}
H.~Coskun, F.~Achilles, R.~DiPietro, N.~Navab, F.~Tombari, Long short-term memory {Kalman} filters: Recurrent neural estimators for pose regularization, in: 2017 IEEE International Conference on Computer Vision (ICCV), October 22-29, Venezia, Italy, 2017, pp. 5525--5533.
\newblock \href {https://doi.org/10.1109/ICCV.2017.589} {\path{doi:10.1109/ICCV.2017.589}}.

\bibitem{art:Liu24}
W.~Liu, Z.~Lai, K.~Bacsa, E.~Chatzi, Neural extended {Kalman} filters for learning and predicting dynamics of structural systems, Structural Health Monitoring 23~(2) (2024) 1037--1052.
\newblock \href {https://doi.org/10.1177/14759217231179912} {\path{doi:10.1177/14759217231179912}}.

\bibitem{art:Olivier21}
A.~Olivier, M.~D. Shields, L.~Graham-Brady, {Bayesian} neural networks for uncertainty quantification in data-driven materials modeling, Computer Methods in Applied Mechanics and Engineering 386 (2021) 114079.
\newblock \href {https://doi.org/10.1016/j.cma.2021.114079} {\path{doi:10.1016/j.cma.2021.114079}}.

\bibitem{art:Cicirello24}
A.~Cicirello, \href{https://arxiv.org/abs/2405.05987}{Physics-enhanced machine learning: a position paper for dynamical systems investigations} (2024).
\newblock \href {http://arxiv.org/abs/2405.05987} {\path{arXiv:2405.05987}}.
\newline\urlprefix\url{https://arxiv.org/abs/2405.05987}

\bibitem{chp:CSAI22}
L.~Rosafalco, A.~Manzoni, S.~Mariani, A.~Corigliano, Combined Model Order Reduction Techniques and Artificial Neural Network for Data Assimilation and Damage Detection in Structures, Springer International Publishing, Cham, 2022, pp. 247--259.
\newblock \href {https://doi.org/10.1007/978-3-030-70787-3\_16} {\path{doi:10.1007/978-3-030-70787-3\_16}}.

\bibitem{art:Seventekidis23}
P.~Seventekidis, D.~Giagopoulos, Model error effects in supervised damage identification of structures with numerically trained classifiers, Mechanical Systems and Signal Processing 184 (2023) 109741.
\newblock \href {https://doi.org/10.1016/j.ymssp.2022.109741} {\path{doi:10.1016/j.ymssp.2022.109741}}.

\bibitem{conti2023reduced}
P.~Conti, G.~Gobat, S.~Fresca, A.~Manzoni, A.~Frangi, Reduced order modeling of parametrized systems through autoencoders and {SINDy} approach: continuation of periodic solutions, Computer Methods in Applied Mechanics and Engineering 411 (2023) 116072.
\newblock \href {https://doi.org/https://doi.org/10.1016/j.cma.2023.116072} {\path{doi:https://doi.org/10.1016/j.cma.2023.116072}}.

\bibitem{art:Opreni23}
A.~Opreni, A.~Vizzaccaro, C.~Touz\'{e}, A.~Frangi, High-order direct parametrisation of invariant manifolds for model order reduction of finite element structures: application to generic forcing terms and parametrically excited systems, Nonlinear Dynamics 111 (2023) 5401--5447.
\newblock \href {https://doi.org/10.1007/s11071-022-07978-3} {\path{doi:10.1007/s11071-022-07978-3}}.

\bibitem{art:Brunton16}
S.~L. Brunton, J.~L. Proctor, J.~N. Kutz, Discovering governing equations from data by sparse identification of nonlinear dynamical systems, Proceedings of the National Academy of Sciences 113~(15) (2016) 3932--3937.
\newblock \href {https://doi.org/10.1073/pnas.1517384113} {\path{doi:10.1073/pnas.1517384113}}.

\bibitem{art:Lathourakis24}
C.~Lathourakis, A.~Cicirello, Physics enhanced sparse identification of dynamical systems with discontinuous nonlinearities, Nonlinear Dynamics 112 (2024) 11237--11264.
\newblock \href {https://doi.org/10.1007/s11071-024-09652-2} {\path{doi:10.1007/s11071-024-09652-2}}.

\bibitem{art:CMAME24}
L.~Rosafalco, P.~Conti, A.~Manzoni, S.~Mariani, A.~Frangi, {EKF–SINDy}: Empowering the extended kalman filter with sparse identification of nonlinear dynamics, Computer Methods in Applied Mechanics and Engineering 431 (2024) 117264.
\newblock \href {https://doi.org/https://doi.org/10.1016/j.cma.2024.117264} {\path{doi:https://doi.org/10.1016/j.cma.2024.117264}}.

\bibitem{art:Dardeno24}
T.~Dardeno, K.~Worden, N.~Dervilis, R.~Mills, L.~Bull, On the hierarchical bayesian modelling of frequency response functions, Mechanical Systems and Signal Processing 208 (2024) 111072.
\newblock \href {https://doi.org/10.1016/j.ymssp.2023.111072} {\path{doi:10.1016/j.ymssp.2023.111072}}.

\bibitem{art:Li24}
S.~Li, Y.~Yang, Data-driven modeling of bifurcation systems by learning the bifurcation parameter generalization, Nonlinear Dynamics 415 (2023) 116204.
\newblock \href {https://doi.org/10.1007/s11071-024-10304-8} {\path{doi:10.1007/s11071-024-10304-8}}.

\bibitem{art:Gotte23}
R.-S. G{\"o}tte, J.~Timmermann, Estimating states and model uncertainties jointly by a sparsity promoting {UKF}, {IFAC}-PapersOnLine 56~(1) (2023) 85--90, 12th IFAC Symposium on Nonlinear Control Systems {NOLCOS} 2022.
\newblock \href {https://doi.org/10.1016/j.ifacol.2023.02.015} {\path{doi:10.1016/j.ifacol.2023.02.015}}.

\bibitem{art:Pal24}
A.~Pal, S.~Nagarajaiah, Sparsity promoting algorithm for identification of nonlinear dynamic system based on unscented kalman filter using novel selective thresholding and penalty-based model selection, Mechanical Systems and Signal Processing 212 (2024) 111301.
\newblock \href {https://doi.org/10.1016/j.ymssp.2024.111301} {\path{doi:10.1016/j.ymssp.2024.111301}}.

\bibitem{art:Julier07}
S.~J. Julier, J.~J. LaViola, On kalman filtering with nonlinear equality constraints, IEEE Transactions on Signal Processing 55~(6) (2007) 2774--2784.
\newblock \href {https://doi.org/10.1109/TSP.2007.893949} {\path{doi:10.1109/TSP.2007.893949}}.

\bibitem{art:Stevens24}
J.~M. Stevens-Haas, Y.~Bhangale, J.~Nathan~Kutz, A.~Aravkin, Learning nonlinear dynamics using kalman smoothing, IEEE Access 12 (2024) 138564--138574.
\newblock \href {https://doi.org/10.1109/ACCESS.2024.3465390} {\path{doi:10.1109/ACCESS.2024.3465390}}.

\bibitem{proc:Wang22}
J.~Wang, J.~Moreira, Y.~Cao, B.~Gopaluni, Time-variant digital twin modeling through the {Kalman}-generalized sparse identification of nonlinear dynamics, in: 2022 American Control Conference (ACC), June 8-10, Atlanta, Georgia, USA, 2022, pp. 5217--5222.
\newblock \href {https://doi.org/10.23919/ACC53348.2022.9867786} {\path{doi:10.23919/ACC53348.2022.9867786}}.

\bibitem{art:Giorgio21}
G.~Gobat, V.~Zega, P.~Fedeli, L.~Guerinoni, C.~Touz{\'e}, A.~Frangi, Reduced order modelling and experimental validation of a {MEMS} gyroscope test-structure exhibiting 1:2 internal resonance, Scientific Reports 11 (2021) 16390.
\newblock \href {https://doi.org/10.1038/s41598-021-95793-y} {\path{doi:10.1038/s41598-021-95793-y}}.

\bibitem{art:Champion19}
K.~P. Champion, S.~L. Brunton, J.~N. Kutz, Discovery of nonlinear multiscale systems: Sampling strategies and embeddings, SIAM Journal on Applied Dynamical Systems 18~(1) (2019) 312--333.
\newblock \href {https://doi.org/10.1137/18M1188227} {\path{doi:10.1137/18M1188227}}.

\bibitem{hastie2015statistical}
T.~Hastie, R.~Tibshirani, M.~Wainwright, Statistical learning with sparsity, Monographs on statistics and applied probability 143~(143) (2015) 8.

\bibitem{art:Lokta1920}
A.~J. Lotka, Analytical note on certain rhythmic relations in organic systems, Proceedings of the National Academy of Sciences 6~(7) (1920) 410--415.
\newblock \href {https://doi.org/10.1073/pnas.6.7.410} {\path{doi:10.1073/pnas.6.7.410}}.

\bibitem{art:Fasel22}
U.~Fasel, J.~N. Kutz, B.~W. Brunton, S.~L. Brunton, Ensemble-{SINDy}: Robust sparse model discovery in the low-data, high-noise limit, with active learning and control, Proceedings of the Royal Society A: Mathematical, Physical and Engineering Sciences 478~(2260) (2022) 20210904.
\newblock \href {https://doi.org/10.1098/rspa.2021.0904} {\path{doi:10.1098/rspa.2021.0904}}.

\bibitem{art:Paolo24}
P.~Conti, J.~Kneifl, A.~Manzoni, A.~Frangi, J.~Fehr, S.~L. Brunton, J.~N. Kutz, \href{https://arxiv.org/abs/2405.20905}{{VENI}, {VINDy}, {VICI}: a variational reduced-order modeling framework with uncertainty quantification} (2024).
\newblock \href {http://arxiv.org/abs/2405.20905} {\path{arXiv:2405.20905}}.
\newline\urlprefix\url{https://arxiv.org/abs/2405.20905}

\bibitem{art:Selkov68}
E.~E. Sel'kov, Self-oscillations in glycolysis 1. a simple kinetic model, European Journal of Biochemistry 4~(1) (1968) 79--86.
\newblock \href {https://doi.org/10.1111/j.1432-1033.1968.tb00175.x} {\path{doi:10.1111/j.1432-1033.1968.tb00175.x}}.

\bibitem{munoz2011introduction}
R.~Mu{\~n}oz-Alicea, Introduction to bifurcations and the hopf bifurcation theorem for planar systems, Operation for {M}athematics 640~(11) (2011).

\bibitem{art:Giorgio23}
G.~Gobat, V.~Zega, P.~Fedeli, C.~Touz\'{e}, A.~Frangi, Frequency combs in a mems resonator featuring 1:2 internal resonance: ab initio reduced order modelling and experimental validation, Nonlinear dynamics 111 (2023) 2991--3017.
\newblock \href {https://doi.org/10.1007/s11071-022-08029-7} {\path{doi:10.1007/s11071-022-08029-7}}.

\bibitem{art:Attilio19}
A.~Frangi, G.~Gobat, Reduced order modelling of the non-linear stiffness in mems resonators, International Journal of Non-Linear Mechanics 116 (2019) 211--218.
\newblock \href {https://doi.org/10.1016/j.ijnonlinmec.2019.07.002} {\path{doi:10.1016/j.ijnonlinmec.2019.07.002}}.

\bibitem{art:MATCONT}
A.~Dhooge, W.~Govaerts, Y.~A. Kuznetsov, {MATCONT}: A {MATLAB} package for numerical bifurcation analysis of {ODEs}, {ACM} Transactions on Mathematical Software 29~(2) (2003) 141–164.
\newblock \href {https://doi.org/10.1145/779359.779362} {\path{doi:10.1145/779359.779362}}.

\bibitem{art:Mirzazadeh18}
R.~Mirzazadeh, A.~Ghisi, S.~Mariani, Statistical investigation of the mechanical and geometrical properties of polysilicon films through on-chip tests, Micromachines 9~(2) (2018).
\newblock \href {https://doi.org/10.3390/mi9020053} {\path{doi:10.3390/mi9020053}}.

\bibitem{art:Brossard20}
M.~Brossard, A.~Barrau, S.~Bonnabel, {AI-IMU} dead-reckoning, IEEE Transactions on Intelligent Vehicles 5~(4) (2020) 585--595.
\newblock \href {https://doi.org/10.1109/TIV.2020.2980758} {\path{doi:10.1109/TIV.2020.2980758}}.

\bibitem{proc:Saccani24}
F.~Saccani, D.~Pau, M.~Amoretti, In-sensor learning for pressure self-calibration, in: 2024 IEEE Sensors Applications Symposium (SAS), 2024, pp. 1--6.
\newblock \href {https://doi.org/10.1109/SAS60918.2024.10636625} {\path{doi:10.1109/SAS60918.2024.10636625}}.

\bibitem{art:Giordano23}
P.~F. Giordano, S.~Quqa, M.~P. Limongelli, The value of monitoring a structural health monitoring system, Structural Safety 100 (2023) 102280.
\newblock \href {https://doi.org/10.1016/j.strusafe.2022.102280} {\path{doi:10.1016/j.strusafe.2022.102280}}.

\bibitem{art:Brunton17}
S.~L. Brunton, B.~W. Brunton, J.~L. Proctor, E.~Kaiser, J.~N. Kutz, Chaos as an intermittently forced linear system, Nature Communications 8~(1) (2017) 19.
\newblock \href {https://doi.org/10.1038/s41467-017-00030-8} {\path{doi:10.1038/s41467-017-00030-8}}.

\bibitem{EKFSINDY_repo}
L.~Rosafalco, P.~Conti, {{EKF-SINDy}}, \url{https://github.com/ContiPaolo/SINDy-EKF} (2024).

\end{thebibliography}

\end{document}